\documentclass[12pt]{article}

\usepackage{amsmath,amsfonts,amssymb}

\newenvironment{demo}{\noindent {\sl Proof}. \ }{\qed \bigskip}

\newtheorem{stheorem}{Theorem }[subsection]
\newtheorem{sdefin}[stheorem]{Definition }
\newtheorem{sprop}[stheorem]{Proposition }
\newtheorem{slemma}[stheorem]{Lemma }

\numberwithin{equation}{section}

\def\half{\frac{1}{2}}
\def\noi{\noindent}
\def\P{\mathbb{P}}
\def\S{\mathbb{S}}
\def\R{\mathbb{R}}
\def\N{\mathbb{N}}

\def\ep{\varepsilon}
\def\ga{\gamma}
\def\omb{\overline \omega}
\def\al{\alpha}
\def\de{\delta}
\def\la{\lambda}
\def\si{\sigma}

\def\cF{{\cal F}}
\def\tix{\tilde x}
\def\esp{{\rm I\mskip-4mu E}}
\def\Hess{{\rm Hess \ }}

\def\qed{\hbox{$\vcenter{\vbox{
   \hrule height 0.4pt\hbox{\vrule width 0.4pt height 6pt
    \kern5pt\vrule width 0.4pt}\hrule height 0.4pt}}$}}

\begin{document}

\title{
\bigskip \bigskip
  Backward Stochastic Differential Equations on Manifolds II}
\author{
\begin{tabular}{c}
Fabrice Blache \\
{\it{Laboratoire de Math\'ematiques}}                 \\
{\it{Universit\'e Blaise Pascal,  63177 Aubi\`ere Cedex, France}}  \\
{\tt{E-mail: Fabrice.Blache@math.univ-bpclermont.fr}}              \\
\end{tabular}}
\date{September 2005}

\maketitle

\vskip .5cm

\centerline{\small{\bf Abstract}}

\begin{center}
\begin{minipage}[c]{330pt}
{\small In \cite{blache03}, we have studied a generalization of the
problem of finding a martingale on a manifold whose terminal value is
known. This article completes the results obtained in the first article
by providing uniqueness and existence theorems in a general framework
(in particular if positive curvatures are allowed), still using 
differential geometry tools.}
\end{minipage}
\end{center}

\bigskip

\section{Introduction}
\label{par11}

\subsection{Setting of the problem}
\label{par11.1}
First, we introduce some notations and definitions.
Unless otherwise stated, we shall work on a fixed finite time interval
$[0;T]$; moreover, $(W_t)_{0 \le t \le T}$ will always denote a Brownian
Motion (BM for short) in $\R^{d_w}$, for a positive integer $d_w$. Moreover,
Einstein's summation convention will be used for repeated indices in
lower and upper position.

Let $(B_t^y)_{0 \le t \le T}$ denote the $\R^d$-valued
diffusion which is the unique strong solution of the following SDE :
\begin{equation}
\left\{
\begin{array}{rcl}
dB_t^y & = & b(B_t^y) dt + \sigma(B_t^y) d W_t  \\
 B_0^y & = & y,
\end{array}
\right.
\label{sde2}
\end{equation}
where $\sigma~: \R^d \rightarrow \R^{d \times d_w}$ and
$b : \R^d \rightarrow \R^d$ are $C^3$ bounded functions with bounded
partial derivatives of order 1, 2 and 3.

Let us recall the problem stated in \cite{blache03}.
We consider a manifold $M$ endowed with a connection $\Gamma$, which
defines an exponential mapping. On $M$, we study the uniqueness and
existence of a solution to the equation (under infinitesimal form)
$$
(M+D)_0
\left\{
      \begin{array}{l}
           X_{t+ d t}= {\rm exp}_{X_t}(Z_t d W_t + f(B_t^y, X_t, Z_t) dt)\\
           X_T=U                         \\
      \end{array}
      \right.
$$
where
$Z_t \in {\cal L}(\R^{d_w}, T_{X_t}M)$ and $f(B_t^y,X_t,Z_t) \in T_{X_t}M$.

When $f=0$, this equation characterizes exactly martingales on $M$ with
terminal value $U$. They are the main tool to solve in a probabilistic way
some PDEs such as the Dirichlet problem or the heat equation. Note in
particular the link with variational problems~: it is
well-known that harmonic mappings (between manifolds) transform Brownian
Motions into martingales (see for instance \cite{kend90}); moreover, these
mappings are critical points of the energy functional and can be used to
model the state of equilibrium of liquid crystals (see the
introduction of \cite{helein} for a brief discussion). \\
In the case of a non-vanishing drift term $f$, the solutions are more
general processes which are linked to more general PDEs and mappings
generalizing harmonic ones. One could see these mappings modelling
the equilibrium state of a liquid crystal in an exterior field equal to
the drift term $f$ in equation $(M+D)_0$. \\
For more details about the links with PDEs, the reader is referred to
Kendall (\cite{kend93}) or Thalmaier (\cite{thal96}); see also the 
introductions of \cite{blache03}, \cite{mpl} and \cite{barmart}.

\bigskip
In local coordinates $(x^i)$, the equation $(M+D)_0$ becomes the following
backward stochastic differential equation (BSDE in short)
$$(M+D)
\left\{
      \begin{array}{l}
          d X_t = Z_t d W_t + \left( - \half \Gamma_{jk}(X_t)
            ([Z_t]^k \vert [Z_t]^j) + f(B_t^y,X_t,Z_t) \right) d t  \\
	  X_T=U.                                               \\
      \end{array}
    \right.   $$
We keep the same notations as in \cite{blache03}~: $( \cdot \vert \cdot)$
is the usual inner product in an Euclidean space, the summation convention
is used, and $[A]^i$ denotes the $i^{th}$ row of any matrix $A$; moreover,
\begin{equation}
\Gamma_{jk} (x) =
\left(
   \begin{array}{c}
       \Gamma^1_{jk}(x)    \\
       \vdots              \\
       \Gamma^n_{jk}(x)
   \end{array}
\right)
\label{christof2}
\end{equation}
is a vector in $\R^n$, whose components are the Christoffel symbols
of the connection. We keep the notations $Z_t$ for a matrix in
$\R^{n \times d_w}$ and $f$ for a mapping from
$\R^d \times \R^n \times \R^{n \times d_w}$ to $\R^n$.
The process $X$ will take its values in a compact set, and a solution of
equation $(M+D)_0$ will be a pair of processes
$(X,Z)$ in $M \times (\R^{d_w} \otimes TM)$ such that $X$ is continuous and
$\esp (\int_0^T \Vert Z_t \Vert_r^2 dt) < \infty$ for a
Riemannian norm $\Vert \cdot \Vert_r$;
in global coordinates $O \subset \R^n$,
$(X,Z) \in O \times \R^{n \times d_w}$ and
$\esp (\int_0^T \Vert Z_t \Vert ^2 dt) < \infty$ (see below for the 
definitions of the norms). 
We gave in \cite{blache03} existence and uniqueness results for the 
solutions of the BSDE
$(M+D)$ in two different frameworks~: firstly when the drift $f$ did not
depend on the process $(Z_t)$, and secondly for a "general" $f$ (i.e.
depending on the process $(Z_t)$), but only
for the Levi-Civita connection and in nonpositive curvatures. \\
In this article we extend these results for "general" drifts $f$ to
other manifolds, without further hypothesis on the curvature as
above. In this
case, we need (unlike in \cite{blache03}) to prove an exponential 
integrability condition like
$$\esp \left( e^{\mu \int_0^T \Vert Z_s \Vert^2 ds} < \infty \right).$$
This allows the construction of a submartingale, which is the crux of the 
matter, but this leads to calculus which is much more intricate than in
\cite{blache03}.
More precisely, we study two cases~: on the one hand,
the case of a general connection (i.e. not depending on the Riemannian
structure defined on the manifold $M$) for which we give only results on
small domains; on the other hand, the case of a manifold endowed with its
Levi-Civita connection and whose sectional curvatures are allowed to be
positive.

In Section \ref{par12}, we prove some estimates and technical lemmas which
are useful in Section \ref{par13}, devoted to prove uniqueness results.
Section \ref{par14} deals with existence results; it is very similar to
Section 4 of \cite{blache03}, so we give the main results generally without
proof, except when the arguments are more complicated than in
\cite{blache03}. To end the article, in Section \ref{par15} we recall 
briefly the links to the Dirichlet problem.

\subsection{Notations and hypothesis}
\label{par11.2}
In all the article, we suppose that a filtered probability space
$(\Omega, {\cal F}, P, ({\cal F}_t)_{0 \le t \le T})$ (verifying the usual
conditions) is given (with $T<\infty$ a deterministic time) on which 
$(W_t)_t$ denotes a $d_w$-dimensional BM. Moreover, we
always deal with a complete Riemannian manifold  $M$ of dimension $n$,
endowed with a linear symmetric (i.e. torsion-free) connection whose
Christoffel symbols $\Gamma^i_{jk}$ are smooth; the connection does not
depend {\it a priori} on the Riemannian structure.

On $M$, $\delta$ denotes the Riemannian distance;
$\vert u \vert_r$ is the Riemannian norm for a tangent vector $u$
and $\vert u' \vert$ the Euclidean norm for a vector $u'$ in $\R^n$.
If $h$ is a smooth real function defined on $M$, its differential is
denoted by $D h$ or $h'$; the Hessian $\Hess h(x)$ is a bilinear
form the value of which is denoted by $\Hess h(x)<u,\overline u>$, for
tangent vectors (at $x$) $u$ and $\overline u$.\\
For $\beta \in \N^*$, we say that a function is $C^\beta$ on a closed set 
$F$ if it is $C^\beta$ on an open set containing $F$.
Recall also that a real function $\chi$ defined on $M$ is said to be
convex if for any $M$-valued geodesic $\gamma$, $\chi \circ \gamma$ is
convex in the usual sense (if $\chi$ is smooth, this is equivalent to
require that $\Hess \chi$ be nonnegative).\\
For a matrix $z$ with $n$ rows and $d$ columns, $^t z$ denotes
its transpose,
$$\Vert z \Vert
= \sqrt{ \hbox{Tr} (z {}^t z)}
= \sqrt{\sum_{i=1}^d \vert [{}^t z]^i \vert ^2}$$
(${\rm Tr}$ is the trace of a square matrix) and
$\Vert z \Vert_r = \sqrt{\sum_{i=1}^d \vert [{}^t z]^i \vert_r ^2}$
where the columns of $z$ are considered as tangent vectors.
The notation $\Psi(x,x') \approx \delta(x,x')^\nu$ means that there is
a constant $c>0$ such that
$$ \forall x,x', \
\frac{1}{c} \ \delta(x,x')^\nu \le \Psi(x,x') \le c \ \delta(x,x')^\nu.$$
Before the general framework, let us give some additional notations
which are specific to the Levi-Civita connection.
In this case, we always assume that the injectivity radius $R$ of $M$
is positive and that its sectional curvatures are bounded above; we let
$K$ be the smallest nonnegative number dominating all the sectional
curvatures, or $0$ if they are all nonpositive.
Finally we recall from \cite{kend90} the definition of a regular geodesic 
ball.
\begin{sdefin}
A closed geodesic ball $\cal B$ of radius $\rho$ and center $p$ is said 
to be regular if
\begin{item}{(i)}
$\rho \sqrt K < \frac{\pi}{2}$
\end{item}

\begin{item}{(ii)}
The cut locus of $p$ does not meet $\cal B$.
\end{item}
\end{sdefin}
For an introductory course in Riemannian
geometry, the reader is referred to Booth\-by (\cite{booth}) and for 
further facts about curvature, to Lee (\cite{lee}).

\bigskip
Let us come back to general connections. Throughout the article, we 
consider an open set $O$ of $M$, relatively
compact in a local chart and an open set $\omega \not= \emptyset$
relatively compact in $O$, verifying that
\medskip
\\ $\bullet$ There is a unique geodesic in $\overline O$, linking any two 
              points of $\overline O$, and depending smoothly on its 
	      endpoints;
\\ $\bullet$ $\omb = \{ \chi \le c \}$, the sublevel set of a smooth
              convex function $\chi$ defined on  $O$.
\medskip    
\noi Note that $O$ will be as well considered as a subset of $\R^n$.
In the case of a general connection, it is well-known that any point $x$
of $M$ has a neighbourhood $O$ for which the first property holds; and
when the Levi-Civita connection is used, the first property is also true
for a regular geodesic ball (see Theorem (1.7) in \cite{kend90}).

\noi Finally we always assume two hypotheses on $f$~:
$$  \exists L>0 , \ \forall b,b' \in \R^d,
  \forall (x,z) \in O \times {\cal L}(\R^{d_w}, T_xM),
  \forall (x',z') \in O \times {\cal L}(\R^{d_w}, T_{x'}M),
$$
\begin{eqnarray}
\left\vert \overset{x'}{\underset{x}{\Vert}} f(b,x,z) - f(b',x',z')
\right\vert_r \le L \Bigg( ( \vert b-b' \vert
& + &\delta(x,x')) (1 + \Vert z \Vert_r + \Vert z' \Vert_r)  \nonumber  \\
& + & \left\Vert \overset{x'}{\underset{x}{\Vert}} z - z' \right\Vert_r
                 \Bigg)
\label{lip11}
\end{eqnarray}
and
\begin{equation}
\exists L_2>0 , \exists x_0 \in O, \forall b \in \R^d,
\vert f(b,x_0,0) \vert_r \le L_2.
\label{upperboundf11}
\end{equation}
The first one is a "geometrical" Lipschitz condition on $f$; this special
form is needed to get an expression which is invariant under changes of
coordinates. The second one means that $f$ is bounded with respect to the 
first argument. Finally we denote by 
$$\overset{x'}{\underset{x}{\Vert}}  z$$ 
the parallel transport (defined by the connection) of the $d_w$ columns
of the matrix $z$ (considered as tangent vectors) along the unique 
geodesic between $x$ and $x'$.

\bigskip
{\bf Remark.} 
In fact, condition \eqref{lip11} can be weakened by 
splitting it into two conditions~:
yet a Lipschitz condition on $b$ and $z$
$$
\left\vert f(b,x,z) - f(b',x,z') \right\vert_r 
\le L \Bigg( \vert b-b' \vert (1 + \Vert z \Vert_r + \Vert z' \Vert_r)  
   + \left\Vert z - z' \right\Vert_r
                                                            \Bigg)
$$
and the following ``monotonicity'' condition on $x$
$$
D \Psi \cdot
     \left(
        \begin{array}{c}
          f(b,x,z)  \\
          f(b,x', \overset{x'}{\underset{x}{\Vert}} z)
        \end{array}
     \right)
 \ge 
 \mu  \Psi(x,x') (1 + \Vert z \Vert) 
$$
for a real constant $\mu$ independent of $b, x, z$ (where $\Psi$ is 
replaced by $\delta$ in the case of the Levi-Civita connection).
This ``monotonicity" condition replaces here the well-known 
monotonicity condition involving the inner product in an Euclidean space (see e.g. Assumption $(4)$ in \cite{darlpard97} or Assumption $(H3)$ in \cite{briandhupard03}). As in these references, we need also 
some additional conditions on $f$, in particular continuity in the $x$ variable.\\
Note that here we have a lower bound on $D \Psi$
(and not an upper bound as in the articles cited above) because in the 
equation $(M+D)$, the drift $f$ is given with a ``plus" sign.

The proof of uniqueness is similar to the one in the Lipschitz case; for the existence result, we approximate $f$ by Lipschitz functions $f_n$ and pass through the limit in equation $(M+D)$; it involves more intricate computations (in particular to keep the assumption that the functions $f_n$ are pointing outward on the boundary of $\omb$, see Subsection \ref{par11.3} below). Details will appear elsewhere.

\bigskip
To end this part, notice that the same letter $C$ will often stand for
different constant numbers.

\subsection{The main result}
\label{par11.3}
Before performing calculations, we give the main theorem of the article.
Let us first introduce a technical but natural hypothesis, which we will
make precise in Section \ref{par14}~:
$$(H) \ \ f \hbox{ is pointing outward on the boundary of } \omb.$$
Then, under the notations above, we can state :
\begin{stheorem}
\label{existence12} 
We consider the BSDE $(M+D)$ with terminal random variable 
$U \in \omb= \{ \chi \le c \}$. We suppose 
that $f$ verifies conditions \eqref{lip11}, \eqref{upperboundf11}, 
and that $\chi$ is strictly convex (i.e. $\Hess \chi$ is positive definite).
\begin{item}
(i) Each point $q$ of $M$ has a 
neighbourhood $O_q \subset O$ such that, if
$\omb \subset O_q$ and $f$ verifies hypothesis $(H)$, then the BSDE 
$(M+D)$ has a unique solution $(X_t,Z_t)_{0 \le t \le T}$ such that $X$ 
remains in $\omb$. The neighbourhood $O_q$ depends on the geometry of the
manifold, but not on the constants $L$ and $L_2$ defined in \eqref{lip11}
and \eqref{upperboundf11}. 
\end{item}
\begin{item}
(ii) If the Levi-Civita connection is used and $\omb \subset {\cal B}$, a
regular geodesic ball, and if $f$ verifies hypothesis $(H)$, then the 
BSDE $(M+D)$ has a unique 
solution $(X_t,Z_t)_{0 \le t \le T}$ with $X$ in $\omb$. Moreover, if 
$\omb={\cal B}$, the hypothesis on the strict convexity of $\chi$ is 
satisfied by taking $\chi = \de^2$.
\end{item}
\end{stheorem}
Theorem \ref{existence12} brings together Theorems \ref{2unicit1} and
\ref{unicit3} and the results of Section \ref{par14}. 
In Theorem \ref{2existence4}, we will extend this theorem to random time
intervals $[0;\tau]$, for stopping times $\tau$ which
verify an exponential integrability condition and "sufficiently small"
drifts $f$.

\bigskip
{\it Acknowledgements : } The author would like to renew his thanks to his
supervisor Jean Picard for his help and his relevant advice, and the 
referees for their suggestions to improve a first version.

\section{Preliminary estimates}
\label{par12}

We first recall elementary results about It\^o's formula and parallel
transport.
Then we give some geometrical
estimates for the distance function on $M \times M$. These results are
nontrivial generalizations (mainly in the case of the Levi-Civita 
connection) of results of \cite{scm} and \cite{mpl}.
In this section, the covariant derivative of a vector field $z(t)$ along
a curve $\gamma_t$ is denoted by $\nabla_{\dot \gamma_t} z(t)$.

\subsection{It\^o's formula on manifolds}

Consider two solutions $(X^1,Z^1)$ and $(X^2,Z^2)$ of equation
$(M+D)$ with terminal values $U^1$ and $U^2$, such that $X^1$
and $X^2$ remain in $O$. Let
$$\tilde X =(X^1,X^2) \ \ \hbox{ and } \ \
    \tilde Z = \left(
                 \begin{array}{c}
                    Z^1   \\
	            Z^2
                  \end{array}
                \right);$$
then, for a smooth function $\Psi$ defined on $O \times O$,
 It\^o's formula is written
	
\begin{eqnarray}
\Psi (\tilde X_t) - \Psi (\tilde X_0)
    &=&  \int_0^t D \Psi (\tilde X_s)
           \left( \tilde Z_s d W_s \right)             \nonumber   \\
    & & +  \int_0^t D \Psi (\tilde X_s)
            \left(
                 \begin{array}{c}
          f(B_s^y,X^1_s,Z^1_s)- \half \Gamma_{jk}(X^1_s)
                                 ([Z^1_s]^k \vert [Z^1_s]^j)  \\
          f(B_s^y,X^2_s,Z^2_s)- \half \Gamma_{jk}(X^2_s)
                                 ([Z^2_s]^{k} \vert [Z^2_s]^{j})
                 \end{array}
            \right)  ds                                \nonumber   \\
    & & +  \half \int_0^t  {\rm Tr}  \left(
                   {}^t \tilde Z_s D^2 \Psi(\tilde X_s) \tilde Z_s
                                          \right) ds    \nonumber  \\
    &=& \int_0^t D \Psi (\tilde X_s)
            \left( \tilde Z_s d W_s \right)            \nonumber   \\
    & & + \half \int_0^t  \left(
      \sum_{i=1}^{d_w} {}^t [{}^t \tilde Z_s]^i \Hess \Psi (\tilde X_s)
                  [{}^t \tilde Z_s]^i \right) ds         \nonumber \\
    & & + \int_0^t  D \Psi(\tilde X_s)
                                   \left(
                                 \begin{array}{c}
                                     f(B_s^y,X^1_s,Z^1_s)  \\
                                     f(B_s^y,X^2_s,Z^2_s)
                                 \end{array}
                                   \right) ds
\label{2ito1}
\end{eqnarray}
where $D^2 \Psi$ is the second order derivatives matrix
(remember also notation \eqref{christof2} and that $[{}^t \tilde Z_s]^i$
denotes the $i^{th}$ column of the
matrix $\tilde Z_s$; it is a vector in $\R^{2n}$).
Moreover, for a smooth function $h$ on $O$ and a solution $(X,Z)$ of
$(M+D)$, we get a similar formula, replacing $\tilde X$ by $X$ and
$\tilde Z$ by $Z$.

\subsection{A comparison result about parallel transports}
The relatively compact set $O$ is considered here as a subset of $\R^n$;
the
following proposition gives in $O$ a comparison result between the parallel
transports defined by the Euclidean structure and the connection.

\begin{sprop}
\label{2proptp}
There is a $C>0$ such that for every $(x,x') \in O \times O$
and  $(z,z') \in T_xM \times T_{x'}M$, we have

$$
\left\vert \overset{x'}{\underset{x}{\Vert}} z - z' \right\vert_r
 \le C
  \left( \vert z-z' \vert + \delta(x,x')(\vert z \vert + \vert z' \vert)
       \right)
$$
and
\begin{equation}
  \vert z-z' \vert \le C
\left( \left\vert \overset{x'}{\underset{x}{\Vert}}
                                    z - z' \right\vert_r
  + \delta(x,x')(\vert z \vert_r + \vert z' \vert_r)
       \right).
 \label{2tp2}
\end{equation}
\end{sprop}

\begin{demo}
The case of the Levi-Civita connection has been treated in \cite{blache03}.
We give the proof here in the case of a general connection, not depending
{\it a priori} on the Riemannian structure.
Using the triangle inequality and the equivalence of the Euclidean and 
Riemannian norms on compact domains, we first remark that it is sufficient 
to prove the existence of $C>0$ such that
\begin{equation}
\forall (x,x') \in O, \forall z \in \R^{n},
\left\vert \overset{x'}{\underset{x}{\Vert}} z - z \right\vert
 \le C \delta(x,x') \vert z \vert.
\label{2tp3}
\end{equation}
In fact, using the linearity property of the parallel transport, it is 
sufficient to prove \eqref{2tp3} for $\vert z \vert =1$.
Define 
\begin{eqnarray*}
\tau : \overline O \times \overline O \times S(0;1)
& \rightarrow & \overline O   \\
(x,x',z)
& \mapsto & \overset{x'}{\underset{x}{\Vert}} z
\end{eqnarray*}
where $S(0;1)$ is the sphere of radius $1$ in $\R^n$.
It is a smooth mapping. Indeed, let $\gamma$ be the geodesic such that 
$\gamma(0)=x$ and $\gamma(1)=x'$; we have supposed that $\gamma$ depends 
smoothly on its endpoints $x$ and $x'$. Moreover, if we note
$$\forall t \in [0;1], \ \
   z(t)=\overset{\gamma_t}{\underset{x}{\Vert}} z,$$
the equation of the parallel transport $\nabla_{\dot \gamma_t} z(t)=0$    
is written in local coordinates
$$ \forall k, \ \ \dot z^k(t) + \Gamma_{jl}^k(\gamma_t) \dot 
    \gamma_t^j z^l(t)=0.$$
The general theory of differential equations gives then the smoothness of 
$\tau$. Now since $(x,x',z)$ is in a compact set, we get, for a constant 
$C$ independent of $x$, $x'$ and~$z$
$$\vert \tau(x,x',z) - \tau(x,x,z) \vert \le C \de(x,x').$$
This is exactly \eqref{2tp3} for $z$ such that $\vert z \vert =1$.
This completes the proof.
\end{demo}

As a consequence, on the relatively compact set $O \subset \R^n$,
\eqref{lip11} becomes

$$  \exists L'>0 , \ \forall b,b' \in \R^d, \
  \forall (x,z), (x',z') \in O \times \R^{n \times d_w},
$$
$$
\vert  f(b,x,z) - f(b',x',z') \vert \le L'
( (\vert b-b' \vert  +  \vert x - x' \vert)
        (1 + \Vert z \Vert + \Vert z' \Vert) +  \Vert z - z' \Vert ).
$$

\subsection{A local lower bound on a Hessian}
\label{par12.2}

In this paragraph, we suppose that there exists a
nonnegative, smooth and convex function $\Psi$ on the product
$\omb \times \omb$ (i.e. convex on an open set containing this set) which
vanishes only on the diagonal $\Delta=\{(x,x)/x \in \omb \}$ ($\omb$ is
said to have $\Gamma$-convex geometry); besides, we
suppose that $\Psi \approx \delta^p$ for a $p \ge 2$. Note that since
$\Psi$ is smooth, $p$ is an even integer (see Remark 2 after the proof of 
Lemma \ref{taylorpsi}.

Let $(x,x')$ be a point in $\omb \times \omb$. For notational
convenience, we keep the same notation $\omb \times \omb$ for the
image of this compact set in the local coordinates considered below, 
and we
write $f$ for $f(b,x,z)$ and $f'$ for $f(b,x',z')$ (note that we use
here the same $b$).
Take a local chart $(\phi,\phi)$ in which $(x,x')$ has coordinates
$(\hat x, \hat x')$; if $({\partial_1}, \ldots, {\partial_{2n}})$
denotes the natural dual basis of these coordinates, then
\begin{equation}
 \left(
   \begin{array}{c}
     f  \\
     f'
   \end{array}
 \right)
 = \sum_{i=1}^n \left( f^i \partial_i + f'^i \partial_{i+n} \right)
\end{equation}
and 
$$
\left( \Hess \Psi (x,x') \right)_{ij} =
 \partial_{ij} \Psi(x,x')
       - \overline \Gamma_{ij}^k(x,x') \partial_k \Psi(x,x') 
 $$
where the $\overline \Gamma_{ij}^k$ are the Christoffel symbols associated
with the product connection. Recall that, if we note
$\overline \Gamma_{ij} =
      (\overline \Gamma_{ij}^1, \ldots, \overline \Gamma_{ij}^{2n})$,
they are given by
$$\overline \Gamma_{ij}(x,x') =
\begin{cases}
(\Gamma_{ij}^1(x), \ldots, \Gamma_{ij}^n(x), 0, \ldots, 0)
& \ \hbox{ if } \  i,j \le n, \\
(0, \ldots, 0, \Gamma_{i-n, j-n}^1(x'), \ldots, \Gamma_{i-n, j-n}^n(x'))
& \ \hbox{ if } \  i,j > n, \\
(0, \ldots 0, 0, \ldots, 0)
& \ \hbox{ if } \  i \le n <j \ \hbox{ or } \  j \le n <  i.
\end{cases}$$
We put
$$
\begin{pmatrix}
 A      & E \\
 {}^t E & B
\end{pmatrix}
_{ij} := \partial_{ij} \Psi  - \overline \Gamma_{ij}^k \partial_k \Psi 
$$
where $A, B$ and $E$ are square matrices of size $n$.

Let $v=(v_1, \ldots, v_{2n})=(\hat x - \hat x', \hat x')$ be new
coordinates and $(\tilde\partial_1, \ldots, \tilde\partial_{2n})$ be the
natural dual basis. Then the diagonal $\Delta$ is given by the equation
$\{ v_1= \cdots = v_n =0 \}$ and for $i=1, \ldots, n$,
$\tilde\partial_i=\partial_i$ and
$\tilde\partial_{i+n} = \partial_i + \partial_{i+n}$.
Therefore in these new coordinates,
$$
 \left(
   \begin{array}{c}
     z  \\
     z'
   \end{array}
 \right)
 = \sum_{i=1}^n \left( (z^i -z'^i) \tilde\partial_i + z'^i
    \tilde\partial_{i+n} \right)
$$
and
$$ {}^t
 \left(
   \begin{array}{c}
     z  \\
     z'
   \end{array}
 \right)
\Hess \Psi(x,x')
 \left(
   \begin{array}{c}
     z  \\
     z'
   \end{array}
 \right)
= {}^t
 \left(
   \begin{array}{c}
     z-z'  \\
     z'
   \end{array}
 \right)
\begin{pmatrix}
 \tilde A      & \tilde E \\
 {}^t \tilde E & \tilde B
\end{pmatrix}
\left(
   \begin{array}{c}
     z-z'  \\
     z'
   \end{array}
 \right),
$$
where the square matrices $\tilde A, \tilde B$ and $\tilde E$ are given
by
$$
\forall i,j = 1, \ldots, n,
\left\{
\begin{array}{l}
\tilde A_{ij} = \tilde \partial_{ij} \Psi
  -\sum_{k=1}^{n} \overline \Gamma_{ij}^k \tilde \partial_k \Psi \\
\tilde B_{ij} = \tilde \partial_{i+n,j+n} \Psi
  -\sum_{k=n+1}^{2n} \overline \Gamma_{i+n,j+n}^k \tilde \partial_k \Psi \\
\tilde E_{ij} = \tilde \partial_{i,j+n} \Psi
  - \sum_{k=1}^{2n} \overline \Gamma_{i,j+n}^k \tilde \partial_k \Psi
      = \tilde \partial_{i,j+n} \Psi.
\end{array}
\right.
$$

\begin{slemma}
\label{taylorpsi}
In $v$-coordinates, we have the following estimates~:
\begin{item}{(i)}
If $i \le n$, $\tilde\partial_i \Psi (x,x') \le C \delta^{p-1}(x,x')$.
\end{item}

\begin{item}{(ii)}
If $i > n$, $\tilde\partial_i \Psi (x,x') \le C \delta^p(x,x')$.
\end{item}
\end{slemma}

\begin{demo}
Let $V$ denote $(x,x')$ in $v$-coordinates and $\hat p(V)$ the projection of
a vector $V$ onto $\Delta$~: if $V=(v_1, \ldots, v_{2n})$, then
$\hat p(V)=(0, \ldots, 0, v_{n+1}, \ldots, v_{2n})$.
Remark first that $\vert V- \hat p(V) \vert \approx \delta(x,x')$;
then use a Taylor expansion of order $p$
of $\Psi$ near the diagonal (remember that $\Psi \approx \delta^p$)~:
\begin{eqnarray*}
\Psi(V) 
& = & \frac{1}{p!} \sum_{1 \le i_1, \dots, i_p \le n} 
        \tilde \partial_{i_1 \dots i_p}
          \Psi (\hat p(V)) (V-\hat p(V))_{i_1} \dots (V-\hat p(V))_{i_p}  \\
& + & \int_0^1\frac{(1-t)^p}{p!} \sum_{1 \le j_1, \dots, j_{p+1} \le n} 
       \tilde \partial_{j_1 \dots j_{p+1}}
        \Psi (V_t) (V- \hat p(V))_{j_1} \dots (V- \hat p(V))_{j_{p+1}} dt
\end{eqnarray*}
where $V_t=\hat p(V)+t(V- \hat p(V))$.
By differentiating this equality, we get the two estimates (note that if 
$i \le n$,  $\tilde \partial_i$ corresponds to the differentiation with
respect to the $i^{th}$ term of $V- \hat p(V)$, and if $i>n$,  
$\tilde \partial_i$ corresponds to the differentiation with
respect to the $i^{th}$ term of $\hat p(V)$).
\end{demo}

{\bf Remark 1 : } By iterating the procedure, one can show in particular
that there is a constant $C>0$ such that for $x$ and $x'$ in 
$\overline O$, $\vert \tilde E_{ij}(x,x') \vert \le C \de(x,x')^{p-1}.$

{\bf Remark 2 : } A consequence of the Taylor expansion in the proof above 
is that $p$ is an even integer. Indeed, it is straightforward that $p$ 
is an integer, and if we change $V$ by $-V+2 \hat p(V)$, we conclude that $p$ 
cannot be odd since $\Psi$ is nonnegative.

\bigskip
We prove now a lower bound on $\Hess \Psi$.
\begin{sprop}
\label{propminorhess}
Suppose that there is a neighbourhood $O_\Delta$ of the diagonal in $\omb$
and a $\eta>0$ such that, for all $x,x'$ in $O_\Delta$ and $z$ in $\R^n$
\begin{equation}
{}^t z \tilde A z \ge \eta \delta^{p-2}(x,x') \vert z \vert^2.
\label{minA}
\end{equation}
Then there are positive constants $\al$ and $\beta$ such that
$$\forall (x,x') \in O_\Delta, \forall (z,z') \in T_xM \times T_{x'}M,$$
\begin{eqnarray}
{}^t
 \left(
   \begin{array}{c}
     z  \\
     z'
   \end{array}
 \right)
\Hess \Psi (x,x')
 \left(
   \begin{array}{c}
     z  \\
     z'
   \end{array}
 \right)
&\ge& \al \delta^{p-2}(x,x') \left\vert \overset{x'}{\underset{x}{\Vert}} 
         z - z' \right\vert_r^2 \nonumber \\
&   & - \beta \Psi(x,x') (\vert z \vert_r^2 + \vert z' \vert_r^2).
\label{minhesspsi}
\end{eqnarray}
\end{sprop}

\begin{demo}
We work in $v$-coordinates. Using the hypothesis on $\tilde A$ and the
nonnegativity of $\tilde B$ (since $\Hess \Psi$ is nonnegative), we write
\begin{eqnarray*}
{}^t
 \left(
   \begin{array}{c}
     z  \\
     z'
   \end{array}
 \right)
\Hess \Psi (x,x')
 \left(
   \begin{array}{c}
     z  \\
     z'
   \end{array}
 \right)
& = &   {}^t (z-z') \tilde A (z-z') + {}^t (z') \tilde B z'        \\
&   &   + {}^t (z-z') \tilde E z'+ {}^t (z') {}^t \tilde E (z-z')             \\                   \\
&\ge& \eta \delta^{p-2}(x,x') \vert z-z' \vert^2
             - 2 \vert {}^t (z-z') \tilde E z' \vert.
\end{eqnarray*}
But using the remark after Lemma \ref{taylorpsi},
$\Vert \tilde E \Vert \le C \delta^{p-1}$, so
\begin{eqnarray*}
\vert {}^t (z-z') \tilde E z' \vert
& \le & C \delta^{p-1}(x,x') \vert z' \vert \ \vert z-z' \vert 
        = C \delta^{p-2}(x,x') \left( \de(x,x') \vert z' \vert \right)
	       \vert z-z' \vert    \\
& \le & \frac{\eta}{2} \delta^{p-2}(x,x') \vert z-z' \vert^2
        + \zeta \Psi(x,x') \vert z' \vert^2
\end{eqnarray*}
where we have used the classical inequality~:
$ ab \le \ep a^2 + (1/\ep) b^2$ for any $\ep>0$, and 
$\Psi \approx \de^p$. This bound gives
$$
{}^t
 \left(
   \begin{array}{c}
     z  \\
     z'
   \end{array}
 \right)
\Hess \Psi (x,x')
 \left(
   \begin{array}{c}
     z  \\
     z'
   \end{array}
 \right)
\ge \frac{\eta}{2} \delta^{p-2}(x,x') \vert z - z' \vert^2 
    - \zeta \Psi(x,x') \vert z' \vert^2.
$$
Using this estimate, the equivalence of norms and Proposition 
\ref{2proptp} gives exactly the result.
\end{demo}

An example of a function $\Psi$ which verifies the hypothesis of
Proposition \ref{propminorhess} is the convex function constructed by
Emery (see Lemma (4.59) in \cite{scm}); in $v$-coordinates as defined
above, it can be written
\begin{equation}
\Psi (v_1, \ldots, v_n, v_{n+1}, \ldots, v_{2n}) =
  \half \left( \ep^2 + \sum_{k=n+1}^{2n} v_k^2 \right)
  \left( \sum_{k=1}^n v_k^2 \right).
\label{defpsiemery}
\end{equation}
Emery has shown that, for $\varepsilon$ sufficiently small, $\Psi$ is 
convex near the diagonal.
Moreover, it also verifies near the diagonal the estimate \eqref{minA}. 
Indeed, with the notations above, we have for $i,j \le n$~:
$$
\tilde A_{ij} (v) =  \tilde \partial_{ij} \Psi (v)
                  -  \sum_{k=1}^{n} \overline \Gamma_{ij}^k(v)
                       \tilde \partial_k \Psi(v).  
$$
An explicit calculation shows that 
$$\vert \tilde \partial_k \Psi (v) \vert 
    \le C \left( \sum_{k=1}^n v_k^2 \right)^\half
  \ \hbox{ and } \ 
\tilde \partial_{ij} \Psi (v) =
\begin{cases}
  (\ep^2 + \sum_{k=n+1}^{2n} v_k^2)
& \hbox{ if } i=j \\   
0 
& \hbox{ otherwise.}   
\end{cases}
$$    
Then, in the sense of matrices, 
$(\tilde \partial_{ij} \Psi)_{ij} \ge \ep^2 I_n$ (where $I_n$ is the 
identity matrix of $\R^n$). But, since
$( \sum_{k=1}^n v_k^2)^\half \approx \vert V-p(V) \vert \approx \de(x,x')$, 
the other matrix 
$$ (\sum_{k=1}^{n} \overline \Gamma_{ij}^k(v) \tilde \partial_k
   \Psi(v))_{ij}$$
is negligible compared to $(\tilde \partial_{ij} \Psi)_{ij}$, near the
diagonal. Thus $(\tilde A_{ij} (v))_{ij} \ge (\ep^2 /2) I_n$ in a
neighbourhood of $\Delta$; this is exactly \eqref{minA} since, in this 
example, $\Psi \approx \de^2$ and $p=2$.

\subsection{Estimates of the derivatives of the distance}
\label{2estimdist}
In this subsection, the Levi-Civita connection is used.
Then the geodesic distance $(x,x') \mapsto \delta(x,x')$
is defined on $M \times M$ and is smooth except on the cut locus and the
diagonal $\{x=x'\}$. We want
to estimate its first and second derivatives when $M \times M$ is endowed
with the product Riemannian metric. If $\tix=(x,x')$ is a point which is
not in the cut locus or the diagonal, there exists a unique minimizing
geodesic $\gamma(t)$, $0 \le t \le 1$, from $x$ to $x'$. If $u_t$ is a
vector of $T_{\gamma(t)}M$, we can decompose $u_t$ as $v_t+w_t$, where 
$v_t$ is the orthogonal projection of $u_t$ on $\dot \gamma(t)$; the 
vectors $v_t$ and $w_t$ are respectively called the tangential and 
orthogonal components of $u_t$. If $u=(u_0,u_1)$ is a vector of 
$T_{\tix} (M \times M)$, $(v_0,v_1)$ and $(w_0,w_1)$ are also called its 
tangential and orthogonal components.
In the sequel, we put
$$y := \sqrt K \frac{\delta(\tix)}{2}.$$
We start with a technical lemma.

\begin{slemma}
\label{maxratio}
For $0<t<\pi$, let $h(t):=\frac{\sin t}{t}$ and
$$
H(t, \beta) := \frac{1- h(t) \cos (t+2\beta)}
                    {\sin^2 \beta + \sin^2(t+\beta)};
$$
then, for $0<y<\frac{\pi}{2}$, we have
$$
\underset{0 \le \beta \le \pi - 2y}{\rm max}(H(2y, \beta))
  = \frac{1+h(2y)}{1+ \cos (2y)}.
$$
\end{slemma}

\begin{demo}
Let $D(t,\beta):= \sin^2 \beta + \sin^2(t+\beta)$ be the denominator of
$H$; then some simple trigonometry gives
\begin{eqnarray*}
D(t,\beta)
&=& \half (1- \cos(2 \beta)) + \half (1- \cos(2t + 2\beta)) \\
&=& 1- \cos t \cdot \cos(t+2 \beta).
\end{eqnarray*}
When we differentiate $H$ with respect to $\beta$, we get
$$ \frac{D(t,\beta) (2h(t) \sin (t+2 \beta))
    -(1- h(t) \cos (t+2\beta))(2 \cos t \cdot \sin (t+ 2 \beta))}
    {D(t,\beta)^2}; 
$$
therefore, we have
$$
\frac{\partial H}{\partial \beta} (t,\beta)
 =  2 (h(t)- \cos t) \frac{\sin(t+2\beta)}{D(t,\beta)^2}.
$$  
In particular, since $h(t) \ge \cos t$ on $[0;\pi]$,
$$sign(\frac{\partial H}{\partial \beta} (t,\beta))
   = sign(\sin(t+2\beta))$$
and the function $H(2y, \cdot)$ reaches its maximum (for
$0 \le \beta \le \pi - 2y$) at $\beta_m = \frac{\pi}{2}-y $. A simple
calculation gives
$$H(2y, \beta_m)= \frac{1+h(2y)}{1+\cos (2y)}.$$
This completes the proof.
\end{demo}

Now we give the estimates on the first two derivatives of the distance.
\begin{slemma}
Let $\tix$ be a point of $M \times M$ such that 
$0< \de(\tix) < \pi / \sqrt K$. Let $u=(u_0,u_1)$ be a vector of 
$T_{\tix}(M \times M)$ and let $v$ and $w$ be its tangential and orthogonal 
components. Then
\begin{equation}
  \vert \delta'(\tix)<u> \vert =
  \left\vert \overset{x'}{\underset{x}{\Vert}} v_0 - v_1 \right\vert _r 
\label{2der1}
\end{equation}
and
\begin{equation}
\Hess \delta (\tix)<u,u>
\ge \frac{1}{\delta(\tix)} \left\vert \overset{x'}{\underset{x}{\Vert}}
        w_0 - w_1  \right\vert_r^2  - \frac{K}{2} \delta(\tix)
  \frac{1+ \frac{\sin(\sqrt K \delta(\tix))}{\sqrt K \delta(\tix)}}
       {1+ \cos(\sqrt K \delta(\tix))}
    \vert w \vert ^2_r.
\label{derkpos}
\end{equation}
\end{slemma}

\begin{demo}
The equality \eqref{2der1} is proved in \cite{blache03}. Now let $J_w(t)$ be
the normal Jacobi field along $\gamma(t)$ satisfying $J_w(0)=w_0$ and
$J_w(1)=w_1$. From (1.1.7) of \cite{mpl} and the proof of Lemma 2.3.1
of \cite{blache03}, we have
\begin{eqnarray*}
 \Hess \delta(\tix)<u,u>
& \ge &  \frac{1}{\delta(\tix)} \int_0^1 \vert \nabla_{\dot \gamma(t)}
         J_w(t) \vert_r^2 dt - K \delta(\tix) \int_0^1
         \vert J_w(t) \vert_r^2 dt                                   \\
& \ge &  \frac{1}{\delta(\tix)}
         \left\vert \overset{x'}{\underset{x}{\Vert}}
                 w_0 - w_1  \right\vert_r^2
         - K \delta(\tix) \int_0^1 \vert J_w(t) \vert_r^2 dt.
\end{eqnarray*}
The problem is to estimate
$\int_0^1 \vert J_w(t) \vert_r^2 dt$. We know from \cite{mpl} that
$$\frac{d^2}{d t^2} \vert J_w(t) \vert_r 
    \ge - K \de^2(\tix) \vert J_w(t) \vert_r $$
at points $t$ such that $J_w(t) \not= 0$. By comparing with the solution
of the corresponding second order differential equation $j$, with 
$j(0)= \vert J_w(0) \vert_r$ and $j(1)= \vert J_w(1) \vert_r $, we have
(recall that $y= \sqrt K \de(\tix)/2$)~:
$$\vert J_w(t) \vert_r \le j(t)= \al \sin (2y t +\beta)   
      \hbox{ on [0;1] },$$
where $\al \ge 0$ and $0 \le \beta \le \pi - 2y$ are defined by
$$\al \sin \beta = \vert w_0 \vert_r, \ \ \
  \al \sin (2y + \beta) = \vert w_1 \vert_r.$$
In particular,
$$\al^2 = \frac{\vert w \vert_r^2}{\sin^2 \beta + \sin^2(2y+\beta)}.$$
Now note that
\begin{eqnarray*}
\int_0^1 \sin^2(2y t +\beta) dt 
& = & \half \left[ t- \frac{\sin (4yt + 2 \beta)}{4y} \right]_0^1 \\
& = & \frac{1}{8y} (4y+ \sin(2 \beta)-\sin(4y+2\beta))            \\
& = & \frac{1}{8y} (4y - 2 \cos(2y+2\beta) \sin(2y))              \\
& = & \half - \half h(2y) \cos(2y+2 \beta). 
\end{eqnarray*}
Then
$$\int_0^1 \vert J_w(t) \vert_r^2 dt
       \le \al^2 \int_0^1 \sin^2(2y t +\beta) dt =
    \frac{1- h(2y) \cos (2y+2\beta)}{\sin^2 \beta + \sin^2(2y+\beta)}
      \cdot \frac{\vert w \vert_r ^2}{2}.$$
According to Lemma \ref{maxratio}, the maximum, for
$0 \le \beta \le \pi - 2y$,
of the first ratio is
$$ \frac{1+h(2y)}{1+\cos (2y)}.$$
Finally,
$$
\Hess \delta (\tix)<u,u>
\ge \frac{1}{\delta(\tix)} \left\vert \overset{x'}{\underset{x}{\Vert}}
        w_0 - w_1  \right\vert_r^2  - \frac{K}{2} \delta(\tix)
  \frac{1+ h(2y)}{1+ \cos(2y)}
    \vert w \vert ^2_r.
$$
This is \eqref{derkpos}.
\end{demo}

For $a$ such that $1<a<2$, we introduce the function $\Psi_a$ defined 
by
\begin{equation}
\Psi_a(x,x')= \sin^a \left( \sqrt K \frac{\delta(x,x')}{2} \right)
            = \sin^a y.
\label{defpsia}
\end{equation}
This is the function we will use to construct the submartingale in the
uniqueness part, when the Levi-Civita connection is used.
In view of It\^o's formula and the uniqueness part, we want to obtain
estimates on $\Hess \Psi_a$.

\begin{slemma}
\label{estimhess}
We have the two following bounds for $\Hess \Psi_a$~:
\begin{item}{(i)}
Let $\beta>1$. There is a constant $\al>0$, and a neighbourhood $V_\beta$ of 
the diagonal, depending on $\al$, such that for any
$\tix \in (V_\beta \setminus \Delta) \cap (\omb \times \omb)$,
\begin{equation}
\Hess \Psi_a (\tix)<u,u> \ge \al \sin^{a-2} (y)
  \left\vert \overset{x'}{\underset{x}{\Vert}} u_0 - u_1  \right\vert_r^2
  - a \beta \frac{K}{2} \Psi_a(\tix) \vert w \vert_r^2.
\label{estimhess1}
\end{equation}
\end{item}

\begin{item}{(ii)}
For any $\tix$ such that $0< \de(\tix)< \pi / \sqrt K$ and any
$u=(u_0,u_1) \in T_x M \times T_{x'} M$,
\begin{equation}
\Hess \Psi_a(\tix)<u,u> \ge -a \frac{K}{2} \Psi_a(\tix) \vert u \vert^2_r.
\label{estimhess2}
\end{equation}
\end{item}
\end{slemma}

\begin{demo}
(i) We can write $\Hess \Psi_a$ explicitly~:
\begin{eqnarray}
\Hess \Psi_a(\tix)<u,u>
& = & \frac{K}{4} a  \sin^{a-2} y \left( (a-1) \cos^2 y - \sin^2 y \right)
            (\delta'(\tix)<u>)^2                                 \nonumber       \\
&   & + a \frac{\sqrt K}{2} \sin^{a-1} y \cos y \ \Hess \delta(\tix)<u,u>.
\label{hesspsi}
\end{eqnarray}
Let us call, in the right part of \eqref{hesspsi}, $T_1$ the first term 
and $T_2$ the second one. We want to bound below these two terms near the
diagonal.
For $T_1$, note that there is a neighbourhood $V_1$ of $\Delta$, in
$\omb \times \omb$, such that, for $\tix \in V_1$ (recall that
$y=\sqrt K \frac{\delta(\tix)}{2}$),
$$ (a-1) \cos^2 y - \sin^2 y \ge \frac{a-1}{2} \ \hbox{ and } \
   \sin^{a-2} y \ge 1;$$
then using \eqref{2der1},
\begin{equation}
T_1 \ge \frac{K}{8} a (a-1) \sin^{a-2} y
   \left\vert \overset{x'}{\underset{x}{\Vert}} v_0 - v_1  \right\vert_r^2.
\label{minhess1}
\end{equation}
For $T_2$, we use \eqref{derkpos} to get the inequality
$$
T_2 \ge
a \frac{K}{4} \sin^{a-2} y \cos y \frac{\sin y}{y}
   \left\vert \overset{x'}{\underset{x}{\Vert}} w_0 - w_1 \right\vert_r^2
-a \frac{K}{2} \Psi_a(\tix) y \cot y \frac{1+ h(2y)}{1+ \cos(2y)}
           \vert w \vert^2_r.
$$
There is again a neighbourhood $V_2$ of $\Delta$, in $\omb \times \omb$,
such that, for $\tix \in V_2$
$$ \cos y \frac{\sin y}{y} \ge \half \ \hbox{ and } \
   y \cot y \frac{1+ h(2y)}{1+ \cos(2y)} \le \beta;$$
therefore, for $\tix \in V_2$,
\begin{equation}
T_2 \ge a \frac{K}{8} \sin^{a-2} y
  \left\vert \overset{x'}{\underset{x}{\Vert}} w_0 - w_1 \right\vert_r^2
- a \beta \frac{K}{2} \Psi_a(\tix) \vert w \vert_r^2.
\label{minhess2}
\end{equation}
Now \eqref{hesspsi} together with \eqref{minhess1} and \eqref{minhess2}
imply the result.

(ii) The equality \eqref{hesspsi} gives
\begin{eqnarray*}
\Hess \Psi_a(\tix)<u,u>
&\ge& -a \frac{K}{4} \Psi_a(\tix) (\delta'(\tix)<u>)^2   \\  
&   & + a \frac{\sqrt K}{2} \sin^{a-1} y \cos y \ \Hess \delta(\tix)<u,u>.
\end{eqnarray*}
Now by \eqref{2der1}, $(\delta'(\tix)<u>)^2 \le 2 \vert v \vert_r^2$;
moreover, we know from estimate (1.1.2) from \cite{mpl} that
$$\Hess \delta(\tix)<u,u> \ge - \sqrt K \tan (y) \vert w \vert_r^2$$
so
$$\Hess \Psi_a(\tix)<u,u> \ge -a \frac{K}{2} \Psi_a(\tix) \vert v \vert_r^2
    -a \frac{K}{2} \Psi_a(\tix) \vert w \vert_r^2.$$
This is \eqref{estimhess2}.
\end{demo}


\section{Uniqueness}
\label{par13}

\subsection{The general method}
Consider two solutions $(X_t,Z_t)_{0 \le t \le T}$ and 
$(X'_t,Z'_t)_{0 \le t \le T}$  of $(M+D)$ such
that $X$ and $X'$ remain in $\omb$ and $X_T=Y_T=U$ (we will say sometimes
"$\omb$-valued solutions of $(M+D)$").
Let
$$\tilde X_s = (X_s,X'_s) \ \ \hbox{ and }
   \tilde Z_s=  \left(
          \begin{array}{c}
             Z_s  \\
             Z'_s
          \end{array}  \right).$$
To prove uniqueness, the idea is to show that the process
$(S_t)_t=(\exp (A_t) \Psi (\tilde X_t))_t$
where 
$$
A_t= \lambda t + \mu \int_0^t (\Vert Z_s \Vert_r^2  +
                           \Vert Z'_s \Vert_r^2 ) ds,
$$
is a submartingale for appropriate nonnegative constants $\lambda$ and
$\mu$, and a suitable function $\Psi$, smooth on $\overline O$. But to
define such a process, we need to consider solutions verifying an
integrability condition. This leads to the following definition.

\begin{sdefin}
If  $\al$ is a positive constant, let $(\cal E_\al)$ be the set of
$\omb$-valued solutions of $(M+D)$ satisfying
$$
  \esp \exp \left( \al \int_0^T \Vert Z_s \Vert_r^2 ds \right)
                 < \infty.
$$
\end{sdefin}
Actually, we will see that for $\al$ small, $(\cal E_\al)$ contains
any solution of the equation $(M+D)$; we will use Lemma 3.4.2 of 
\cite{blache03}
which we recall here and which is valid for a general connection.

\begin{slemma}
\label{2Emu}
Suppose that we are given a positive
constant $\al$ and a $C^2$ function $\phi$ on $\omb$ satisfying
$C_{min} \le \phi(x) \le C_{max}$ for some positive $C_{min}$ and
$C_{max}$. Suppose moreover that
$\Hess \phi + 2 \al \phi \le 0$ on $\omb$; this means that
\begin{equation}
  \Hess \phi(x)<u,u>+2 \al \phi(x) \vert u \vert_r^2 \le 0.
\label{2min1}
\end{equation}
Then, for every $\ep>0$, any $\omb$-valued solution of $(M+D)$ belongs
to $(\cal E_{\al - \ep})$.
\end{slemma}
Once we get the integrability property, it remains to show that $(S_t)_t$ 
is indeed a submartingale. We now turn to that problem. An application of 
It\^o's formula \eqref{2ito1} gives
\begin{eqnarray}
S_t - S_0
& = &\int_0^t e^{A_s} d (\Psi (\tilde X_s))
+ \int_0^t e^{A_s} (\lambda + \mu (\Vert Z_s \Vert_r^2
         + \Vert Z'_s \Vert_r^2)) \Psi (\tilde X_s) ds    \nonumber \\
& = &\int_0^t e^{A_s} D \Psi (\tilde X_s)
          \left( \tilde Z_s d W_s \right)                 \nonumber \\
&   & + \half \int_0^t e^{A_s}
          \left(  \sum_{i=1}^{d_w} {}^t [{}^t \tilde Z_s]^i
   \Hess \Psi (\tilde X_s) [{}^t \tilde Z_s]^i \right) ds \nonumber \\
&  & + \int_0^t e^{A_s} D \Psi(\tilde X_s)
                                   \left(
                                 \begin{array}{c}
                                     f(B_s^y,X_s,Z_s)  \\
                                     f(B_s^y,X'_s,Z'_s)
                                  \end{array}
                                   \right) ds             \nonumber \\
&  & + \int_0^t e^{A_s} \Psi (\tilde X_s) (\lambda
     + \mu (\Vert Z_s \Vert_r^2 + \Vert Z'_s \Vert_r^2)) ds.
\label{2ito3}
\end{eqnarray}
It is clear that the submartingale property will hold if we show the
nonnegativity of the sum
\begin{eqnarray}
\half \sum_{i=1}^{d_w} {}^t [{}^t \tilde Z_t]^i
        \Hess \Psi (\tilde X_t) [{}^t \tilde Z_t]^i
& + &  D \Psi(\tilde X_t)
                                   \left(
                                 \begin{array}{c}
                                     f(B_t^y,X_t,Z_t)  \\
                                     f(B_t^y,X'_t,Z'_t)
                                  \end{array}
                                   \right)             \nonumber     \\
& + & (\lambda + \mu (\Vert Z_t \Vert_r^2 + \Vert Z'_t \Vert_r^2))
            \Psi(\tilde X_t).
\label{2pos}
\end{eqnarray}
Before giving the details of the proof, we first recall the following upper 
bound on the second term in \eqref{2pos}. 
\begin{slemma}
\label{2upbddpsi}
Suppose that $\Psi$ is a smooth nonnegative function on $\omb \times \omb$,
vanishing only on the diagonal and such that $\Psi \approx \delta^p$ (we 
have seen then that $p$ is an even positive integer).
Then, in local coordinates, there is $C>0$ such that,
for all $x,x'$ in $\omb$, $z,z'$ in $\R^{n \times d_w}$ and $b \in \R^d$
\begin{eqnarray}
 \left\vert D \Psi \cdot
                                   \left(
                                 \begin{array}{c}
                                     f(b,x,z)  \\
                                     f(b,x',z')
                                  \end{array}
                                   \right)
 \right\vert
& \le &
 C  \delta^{p-1}(x,x') \left( \delta(x,x')
(1 + \Vert z \Vert + \Vert z' \Vert) + \Vert z - z' \Vert \right)
                                                            \nonumber \\
& \le &
 C_\ep \Psi(x,x') (1 + \Vert z \Vert_r + \Vert z' \Vert_r)  \nonumber \\
&     & 
+ \ep  \delta^{p-2}(x,x') \left\Vert 
     \overset{x'}{\underset{x}{\Vert}} z - z' \right\Vert_r^2.
\label{2majdpsi}
\end{eqnarray}
\end{slemma}

\begin{demo}
The first inequality is a consequence of Lemma 3.2.1 in
\cite{blache03}. The second results from the algebraic inequality
$$ \de(x,x') \Vert z-z' \Vert \le \frac{1}{\ep} \de^2(x,x') 
      + \ep \Vert z-z' \Vert^2,$$
the hypothesis $\Psi \approx \de^p$, the equivalence of Riemannian and
Euclidean norms, and the inequality \eqref{2tp2}. 
\end{demo}

\subsection{The case of a general connection}
In this paragraph, the connection does not depend {\it a priori} on the
Riemannian structure.

\begin{slemma}
Let $\mu_0$ be a positive real number.
Each point $q$ of the manifold $M$ has a neighbourhood $O^{int}_q$,
relatively compact in $M$ and depending on the geometry of the manifold
(but not on $f$), with the following property~: if $(X,Z)$ is a
solution of equation $(M+D)$ such that $X \in \overline{O^{int}_q}$, then
it belongs to $({\cal E}_{\mu_0})$.
\label{2integ}
\end{slemma}

\begin{demo}
We want to apply the result of Lemma \ref{2Emu}; we consider the function
$g(x)= \cos(a \vert x \vert)$ in a normal system of coordinates centered
at $q$, where $a>0$ and $x$ is such that
$\vert x \vert \le r_0 < \pi / 2a$ ($a$ and $r_0$ are to be determined).
Firstly, $g$ is smooth and bounded, below and above, by positive
constants.\\
Now, for $\mu>\mu_0$, we prove that we can find $a$ large enough such that
$g$ verifies the inequality \eqref{2min1} with $\al=\mu$.
Using the equivalence of the Riemannian and the Euclidean norms on the 
domain on which we are working (by shrinking it if necessary), we know
that there is a $C_r>0$ such that 
$\vert \cdot \vert_r^2 \le C_r \vert \cdot \vert^2$. In particular, it 
is sufficient to show that for any $u \in \R^n$,
\begin{equation}
\Hess g(x)<u,u> + 2 C_r \mu g(x) \vert u \vert^2 \le 0.
\label{minor1}
\end{equation}
Recall that, for $u=(u_1, \ldots, u_n)$
$$\Hess g(x)<u,u>= \sum_{i,j} \left( D_{ij}g(x) - \Gamma_{ij}^k(x) D_kg(x)
      \right) u_i u_j.$$
It is easy to compute the partial derivatives of $g$
\begin{eqnarray*}
D_k g(x)     &=& -a \sin(a \vert x \vert) \frac{x_k}{\vert x \vert} \\
D_{ij} g(x)  &=& -a^2 \cos(a \vert x \vert) \frac{x_i x_j}{\vert x \vert^2}
                 -a \sin(a \vert x \vert) \frac{\delta_{ij}}{\vert x \vert}
                 +a \sin(a \vert x \vert) \frac{x_i x_j}{\vert x \vert^3}
\end{eqnarray*}
(where $\delta_{ij}=1$ if $i=j$, $0$ otherwise).
So \eqref{minor1} is equivalent to showing that the following matrix
\begin{equation}
\left( -a^2 g(x) \frac{x_i x_j}{\vert x \vert^2}
   +a \sin(a \vert x \vert) \frac{x_i x_j}{\vert x \vert^3}
  + a \frac{\sin(a \vert x \vert)}{\vert x \vert}
       \left( \Gamma_{ij}^k(x) x_k - \delta_{ij} \right)
  + 2 C_r \mu g(x) \delta_{ij} \right)_{ij}
\label{matrix1}
\end{equation}
is nonpositive (as a symmetric matrix).\\
We put $a=\sqrt{2 \pi C_r \mu}$ and 
$k(x)=a \sin(a \vert x \vert) / \vert x \vert$; note that 
$k(x) \rightarrow a^2 >0$ if $x \rightarrow 0$. Then
$$
 - \half k(x) + 2 C_r \mu g(x) 
= a^2 \left( - \frac{\sin(a \vert x \vert)}{2 a \vert x \vert}
    + \frac{1}{\pi} g(x) \right) \le 0,
$$
since $0 \le a \vert x \vert < \pi/2$ and $\sin t \ge 2t / \pi$ 
for $t \in [0; \pi/2]$. 
Therefore the matrix 
\begin{equation}
\left( -\half k(x) + 2 C_r \mu g(x) \right) (\de_{ij})_{ij}
\label{matrixneg1}
\end{equation}
is nonpositive. Moreover we have 
\begin{eqnarray*}
 \left( -a^2 g(x) \frac{x_i x_j}{\vert x \vert^2} 
           +a \sin(a \vert x \vert) \frac{x_i x_j}{\vert x \vert^3}
                                               \right)_{ij} 
& = & a^4 \left( \frac{-a \vert x \vert g(x) + \sin(a \vert x \vert)}
                       {a^3 \vert x \vert^3} \right) (x_i x_j)_{ij} \\
&\le& a^4 (x_i x_j)_{ij}
\end{eqnarray*}
if $\vert x \vert \le r_1$, for $(-t \cos t + \sin t)/t^3 \le 1$ in a
neighbourhood of $0$ and the matrix $(x_i x_j)_{ij}$ is easily seen to be
nonnegative. Thus the matrix
\begin{equation}
 \left( -a^2 g(x) \frac{x_i x_j}{\vert x \vert^2} 
           +a \sin(a \vert x \vert) \frac{x_i x_j}{\vert x \vert^3}
               - \frac{1}{4} k(x) \de_{ij} \right)_{ij} 
\label{matrixneg2}
\end{equation}
is nonpositive if $\vert x \vert \le r_2 \le r_1$. Now, for 
$\vert x \vert \le r_3$, the matrix
$$
 \left( a \frac{\sin(a \vert x \vert)}{\vert x \vert}
       \left( \Gamma_{ij}^k(x) x_k - \delta_{ij} \right)
    \right)_{ij}
$$
is less than $- (3/4) k(x) (\delta_{ij})_{ij}$. Using this and the 
nonpositivity of the matrices \eqref{matrixneg1} and \eqref{matrixneg2}, we
conclude that, for $\vert x \vert \le r_0= \inf(r_2,r_3, \pi / (4a))$, the 
matrix \eqref{matrix1} is nonpositive. 
 
This shows the inequality \eqref{minor1}.
As a consequence of Lemma \ref{2Emu}, if we take 
 $\overline{O^{int}_q}= \exp_q(\overline{B(0,r_0)})$, any solution $(X,Z)$
of the equation $(M+D)$, with $X \in \overline{O^{int}_q}$, is in 
$(\cal E_{\mu - \ep})$ for any $\ep>0$; it is in particular in
$({\cal E}_{\mu_0})$.\\
Remark that the drift $f$ interferes in the proof (given in \cite{blache03})
of Lemma \ref{2Emu}, but the result of this lemma is independent of 
$f$. Therefore the neighbourhood $\overline{O^{int}_q}$ does not depend 
on $f$.
\end{demo}

The following theorem gives the uniqueness property in the case of a
general connection, independent of the Riemannian structure.
\begin{stheorem}
Each point $q$ of the manifold $M$ has a neighbourhood $O_q$,
relatively compact in $M$, and depending only on the geometry of the 
manifold (but not on the constants $L$ and $L_2$ in the conditions
\eqref{lip11} and \eqref{upperboundf11} on $f$), with the following property~:
for a given terminal value $U$ in $O_q$, there is at most one
solution $(X,Z)$ to the equation $(M+D)$ such that $X$ remains in $O_q$.
\label{2unicit1}
\end{stheorem}

\begin{demo}
Fix $q \in M$; according to Lemma (4.59) of \cite{scm}, $q$ has a
neighbourhood $O_1$, relatively compact in $M$, with $\Gamma$-convex
geometry, with the convex function $\Psi$ defined at the end of Subsection
\ref{par12.2}. We recall that this function $\Psi$ verifies the estimate
\eqref{minhesspsi} near the diagonal (shrinking $O_1$ if necessary, we
can suppose that it holds on $O_1 \times O_1$) and $\Psi \approx \de^2$.
Taking $p=2$ in \eqref{2majdpsi} and $\ep=\al$ as in \eqref{minhesspsi},
we get
\begin{eqnarray*}
U_t
& = & \half \sum_{i=1}^{d_w} {}^t [{}^t \tilde Z_t]^i
        \Hess \Psi (\tilde X_t) [{}^t \tilde Z_t]^i
              +  D \Psi(\tilde X_t)
                                   \left(
                                 \begin{array}{c}
                                     f(B_t^y,X_t,Z_t)  \\
                                     f(B_t^y,X'_t,Z'_t)
                                  \end{array}
                                   \right)      \nonumber \\
& \ge &
    - \frac{\beta}{2} \Psi(\tilde X_t) 
             (\Vert Z_t \Vert_r^2 + \Vert Z'_t \Vert_r^2)
-C_\al (1 + \Vert Z_t \Vert_r + \Vert Z'_t \Vert_r) \Psi(\tilde X_t).
\end{eqnarray*}
Thus for $\lambda_0$ and $\mu_0$ large enough,
$$ U_t + (\lambda_0 + \mu_0 (\Vert Z_t \Vert_r^2 + \Vert Z'_t \Vert_r^2)
    \Psi(\tilde X_t) \ge 0. 
$$
Note that it suffices to take $\mu_0 > \beta / 2$, and that $\beta$ depends
only on $\Psi$. Therefore $\mu_0$ does not depend on the drift $f$.
For these $\lambda_0$ and $\mu_0$, it results from \eqref{2ito3} that 
the process $(S_t)_t = (\exp (A_t) \Psi (\tilde X_t))_t$
is a submartingale under the additional integrability condition
$$ \esp \left( e^{\mu_0 \int_0^t
   (\Vert Z_s \Vert_r^2+ \Vert Z'_s \Vert_r^2 ) ds} \right) < \infty.$$
According to Lemma \ref{2integ}, this inequality holds if $X$ remains in
a small compact neighbourhood $\overline{O^{int}_q}$ of $q$.
Thus the process $(S_t)_t$ is indeed a submartingale if $X$ remains
in $\overline{O_1 \cap O^{int}_q}$.
Now the conclusion is classical~: since the submartingale $S$ is 
nonnegative and has terminal
value $0$, it vanishes identically. Therefore $\Psi(\tilde X_t)=0$ for 
any  $t$, and the definition of $\Psi$ leads to $X_t=X'_t$. The
proof is completed since we consider continuous processes.
\end{demo}

\subsection{Uniqueness for regular geodesic balls}
In this subsection, the connection used is Levi-Civita's one.
The way to get uniqueness is similar as the one in the preceding
subsection; the difference is that we can obtain more precise results
by achieving explicit calculations.

An accurate examination of formula \eqref{derkpos} together with
\eqref{2pos} shows that the $\mu$ (in \eqref{2pos}) needed to balance, near
the diagonal, the second term in the Hessian is approximately
$\frac{K}{4}$; this implies that the processes should be in
$({\cal E}_{\frac{K}{4}})$. As the next lemma shows, this condition
requires that we work on smaller domains than in the nonpositive-curvature
case treated in \cite{blache03}.
Let $o \in M$ and ${\cal B}_\rho$ denote the regular geodesic ball of radius
$\rho$ and centered at $o$.

\begin{slemma}
\label{leminteg}
Suppose that $0<\gamma<1$ is such that 
$\rho = \gamma \frac{\pi}{2 \sqrt K}$; if we are given two solutions 
$X, X'$ of $(M+D)$ remaining in ${\cal B}_\rho$, then for any 
$0<e<1/\gamma$,
\begin{item}{(i)}
$X$ and $X'$ belong to $({\cal E}_{e \frac{K}{2}})$;
\end{item}

\begin{item}{(ii)}
$\esp \exp \left( e \frac{K}{4} \int_0^T (\Vert Z_s \Vert_r^2
            + \Vert Z'_s \Vert_r^2) ds \right) < \infty.$
\end{item}
\end{slemma}

\begin{demo}
(i)
Let $\varphi(x)=\cos(\beta \sqrt K \delta(o,x))$ with
$1 / \gamma > \beta >1$.
Then on ${\cal B}_\rho$, $\varphi$ is
smooth and there is $c>0$ such that $1 \ge \varphi(x) \ge c >0.$
Moreover,
\begin{eqnarray*}
\Hess \varphi(x)<u,u>
& = & -K \beta^2 \varphi(x) \left( \delta'(o,x)<u> \right)^2 \\
&   & -\sqrt K \beta \sin( \sqrt K \beta \delta(o,x)) \cdot
                 \Hess \delta (o,x) <u,u>                     \\
&\le& - K \beta^2 \varphi(x) \vert v \vert^2_r                \\
&   & - K \beta \sin( \sqrt K \beta \delta(o,x))
        \cot (\sqrt K \delta(o,x)) \vert w \vert^2_r. 
\end{eqnarray*}
The last inequality is a consequence of \eqref{2der1} and the estimation
$(1.2.2)$ of \cite{mpl}~:
\begin{equation}
\Hess \delta (o,x) <u,u> \ge \sqrt K \cot (\sqrt K \delta(o,x)) 
\vert w \vert^2_r.
\label{minhessdelta}
\end{equation}
As the function cotangent is decreasing on
$]0;\frac{\pi}{2}]$ and $\beta>1$, finally we get
$$
\Hess \varphi(x)<u,u> \le -K \beta \varphi(x) \vert u \vert^2_r.
$$
Now we can apply Lemma \ref{2Emu}~: $X$ and $X'$ belong to 
$({\cal E}_{K \beta /2} - \ep)$ for any $\ep>0$; but since $\beta$
is arbitrary between $1$ and $1/ \ga$, $X$ and $X'$ belong to
$({\cal E}_{e K / 2})$ for $0<e<1 / \gamma$.

(ii) It is an immediate consequence of $(i)$ and Cauchy-Schwarz 
inequality.
\end {demo}

From now on in this subsection, we consider only processes remaining in
$\omb={\cal B}_\rho$ with $\rho= \gamma \frac{\pi}{2 \sqrt K}$ 
($0<\gamma<1$) fixed. According to Lemma \ref{leminteg}, fix 
$1<e<1 / \gamma$ such that
$$\esp \exp \left( e \frac{K}{4} \int_0^T (\Vert Z_s \Vert_r^2
            + \Vert Z'_s \Vert_r^2) ds \right) < \infty;$$
then fix $a$ such that $1<a<1+ (e-1) / 2$.
We are going to prove now that the sum \eqref{2pos} is nonnegative for
$\Psi = \Psi_a$, the function defined in \eqref{defpsia}, and for
suitable constants $\la$ and $\mu$.

\begin{sprop}
Let $\Psi(\tix)= \sin^a \left( \sqrt K \delta(\tix) / 2 \right)$ and
$\si, \tau$ stopping times such that $0 \le \si < \tau \le T$. Then, if
$\P [ \forall t \in [\si; \tau], X_t \notin \Delta ] = 1$, the sum
\eqref{2pos} is nonnegative for $\si \le t \le \tau$ and
$(S_t)_t = (\exp(A_t) \Psi(\tilde X_t))_t$
is a submartingale on the random interval $[\si; \tau]$, for $\la$ large 
enough and $\mu = e K/4$ with $e$ as above ($\la$ and $\mu$ don't 
depend on $\si$ nor $\tau$).
\end{sprop}

\begin{demo}
First remark that $\Psi$ is smooth out of the diagonal but not on the
diagonal, so we will consider in this proof only $\tix =(x,x') \notin 
\Delta$. Moreover, we put again $y:=\sqrt K \delta(\tix) /2$ and
$h(y):= \sin y / y$.
We complete the proof in two steps~: first we work near the diagonal and
then away from it.

Near the diagonal, we use part (i) of Lemma \ref{estimhess} with
$\beta>1$ such that $a \beta = 1+ (e-1)/2$ (this is possible
since $a<1+(e-1)/2$).
Let $b \in \R^d$, $x, x' \in (\overline V_\beta \setminus \Delta) \cap 
                      ({\cal B}_\rho \times {\cal B}_\rho)$,
$z \in T_xM$, $z' \in T_{x'}M$ and $u =(z,z')$. On the one hand, we have
from \eqref{estimhess1}
$$
\Hess \Psi (\tix)<u,u> \ge \al \sin^{a-2} (y)
  \left\vert \overset{x'}{\underset{x}{\Vert}} z - z'  \right\vert_r^2
  - \left( 1+\frac{e-1}{2} \right) \frac{K}{2} \Psi(\tix) (\vert z \vert_r^2
         + \vert z' \vert_r^2).
$$
On the other hand, for $z \in {\cal L}(\R^{d_w}, T_xM)$ and
$z' \in {\cal L}(\R^{d_w}, T_{x'}M)$,
\begin{eqnarray*}
  \left\vert D \Psi (\tix) \cdot     \left(
                                 \begin{array}{c}
                                     f(b,x,z) \\
                                     f(b,x',z')
                                  \end{array}
                                     \right)  \right\vert
& = & a \sin^{a-1}(y) \cos y \cdot \frac{\sqrt K}{2}
   \left\vert \delta'(\tix) \cdot              \left(
                                 \begin{array}{c}
                                     f(b,x,z)  \\
                                     f(b,x',z')
                                  \end{array}
                                    \right) \right\vert                 \\
&\le& a \sin^{a-1}(y) \frac{\sqrt K}{2}                         \\
&   & \hskip 1cm \times L \left( \delta(\tix) 
            (1+\Vert z \Vert_r + \Vert z' \Vert_r )
      + \left\Vert \overset{x'}{\underset{x}{\Vert}} z - z'  \right\Vert_r
         \right)                                                        \\
&\le& a \frac{\pi}{2} \Psi(\tix) L(1+ \Vert z \Vert_r + \Vert z' \Vert_r )
      + C_1 \sin^{a-1}y  \left\Vert \overset{x'}{\underset{x}{\Vert}} 
             z - z'  \right\Vert_r \\
&\le& C \Psi(\tix) + \frac{e-1}{2} \frac{K}{4} \Psi(\tix) 
           (\Vert z \Vert^2_r + \Vert z' \Vert^2_r) \\
&   & + \frac{\al}{2} \sin^{a-2} y
       \left\Vert \overset{x'}{\underset{x}{\Vert}} z - z'  \right\Vert_r^2. 
\end{eqnarray*}
The first inequality above is a consequence of \eqref{2der1} and \eqref{lip11}, 
and the second one of the inequality $y \le (\pi / 2) \sin y$, because 
$0<y<\pi /2$. The last one uses classical inequalities.
Therefore, summing the terms 
$\Hess \Psi (\tilde X_t) <[{}^t \tilde Z_t]^i, [{}^t \tilde Z_t]^i>$ for
$i = 1, \ldots, d_w$, we get 
\begin{eqnarray*}
\half
{\rm Tr} \left( ^t \tilde Z_t \Hess
           \Psi (\tilde X_t) \tilde Z_t \right) & &         \\
 +   D \Psi(\tilde X_t)
                                   \left(
                                 \begin{array}{c}
                                     f(B_t^y,X_t,Z_t)  \\
                                     f(B_t^y,X'_t,Z'_t)
                                  \end{array}
                                   \right)
&\ge& - C \Psi(\tilde X_t) - e \frac{K}{4} \Psi(\tilde X_t) 
      (\Vert Z_t \Vert_r^2  + \Vert Z'_t \Vert_r^2).
\end{eqnarray*}
Taking $\lambda \ge C$ and $\mu=e \frac{K}{4}$, the sum \eqref{2pos} is 
nonnegative on 
$(\overline V_\beta \setminus \Delta) \cap ({\cal B}_\rho 
   \times {\cal B}_\rho)$.

It remains to show the nonnegativity of \eqref{2pos} on
$({\cal B}_\rho \times {\cal B}_\rho) \setminus \overline V_\beta$. Using 
here  estimation \eqref{estimhess2}, we get, since $a<e$ (using the 
definition of $a$ and $e$)  
\begin{eqnarray*}
\half
{\rm Tr} \left( ^t \tilde Z_t \Hess
                               \Psi (\tilde X_t) \tilde Z_t \right)  
& + & (\lambda + e \frac{K}{4} (\Vert Z_t \Vert_r^2 + \Vert Z'_t \Vert_r^2))
            \Psi(\tilde X_t) \\
& \ge& 
\lambda \Psi(\tilde X_t) + (e-a) \frac{K}{4} \Psi(\tilde X_t) 
	(\Vert Z_t \Vert_r^2  + \Vert Z'_t \Vert_r^2).   
\end{eqnarray*}
But we also have, for $x, x'$ in 
$({\cal B}_\rho \times {\cal B}_\rho) \setminus \overline V_\beta$,
\begin{eqnarray*}
  \left\vert D \Psi(\tix) \cdot     \left(
                                 \begin{array}{c}
                                     f(b,x,z) \\
                                     f(b,x',z')
                                  \end{array}
                     \right) \right\vert
& \le & C_2 (1+ \Vert z \Vert_r + \Vert z' \Vert_r)                \\
& \le & C_3 + (e-a) \left(\min_{({\cal B}_\rho \times {\cal B}_\rho) 
                                \setminus \overline V_\beta}
 \Psi \right) \frac{K}{4} (\Vert z \Vert^2_r + \Vert z' \Vert^2_r).
\end{eqnarray*}
Then taking $\lambda \ge C_3 / (\min \Psi)$ gives the nonnegativity 
of the sum \eqref{2pos} outside $\overline V_\beta$.
Then the sum \eqref{2pos} is always nonnegative, and the submartingale
property follows.
\end{demo}

The next proposition extends the preceding result to processes that are 
allowed to live in the diagonal.  
\begin{sprop}
If $(\tilde X_t)_t$ remains in 
$\overline {\cal B}_\rho \times \overline {\cal B}_\rho$, then the process 
$(S_t)_{0 \le t \le T}$ is a submartingale.
\end{sprop}

\begin{demo}
Let $S_t^\ep:= S_t \vee \ep$ for $\ep>0$; we begin by proving that the 
process $(S_t^\ep)_t$ is a submartingale.
Indeed, let 
$$\tau'_0:=0, \ 
  \tau_k:=\inf \{ t \ge \tau'_{k-1} : S_t \le \frac{\ep}{2} \} \ 
      \hbox{ and } \ \tau'_k:=\inf \{ t \ge \tau_k : S_t \ge \ep \}.$$
Then it is sufficient to 
show that for $k$ fixed, $(S_t^\ep)$ is a submartingale on the random time 
intervals $[\tau_k;\tau'_k]$ and $[\tau'_k;\tau_{k+1}]$.
Let $\si$ and $\tau$ be stopping times; if  
$\tau_k \le \si < \tau \le \tau'_k$ and $A \in {\cal F}_\si$ then 
$\esp \left( S_\tau^\ep 1_A \right) = \esp \left( S_\si^\ep 1_A \right)$
since $S_\si^\ep=S_\tau^\ep=\ep$.
Now suppose that $\tau'_k \le \si < \tau \le \tau_{k+1}$ and 
$A \in {\cal F}_\si$; then 
$\esp \left( S_\tau 1_A \right) \ge \esp \left( S_\si 1_A \right)$ 
since $(S_u)$ is a submartingale out of the diagonal.
In particular,
$$\esp \left( S_\tau^\ep 1_A 1_{ \{S_\si>\ep \} } \right) \ge
  \esp \left( S_\tau 1_A 1_{ \{S_\si>\ep \} } \right) \ge
  \esp \left( S_\si 1_A 1_{ \{S_\si >\ep \} } \right) =
  \esp \left( S_\si^\ep 1_A 1_{ \{S_\si>\ep \} } \right);$$
moreover,
$$\esp \left( S_\tau^\ep 1_A 1_{ \{S_\si \le \ep \} } \right) \ge
  \esp \left( S_\si^\ep 1_A 1_{ \{S_\si \le \ep \} } \right).$$
So 
$\esp \left( S_\tau^\ep 1_A \right) \ge \esp \left( S_\si^\ep 1_A \right)$.
Therefore, $(S_t^\ep)$ is indeed a submartingale.
To conclude, it suffices to apply Lebesgue's dominated convergence theorem
with $\ep \rightarrow 0$, since $\sup_t S_t$ is integrable.
\end{demo}

The uniqueness property now follows as usual.

\begin{stheorem}
\label{unicit3}
Suppose that the drift $f$
verifies conditions \eqref{lip11} and \eqref{upperboundf11}, and that
$\omb = {\cal B}_\rho$ is a regular geodesic ball. Then,
for a given terminal value $U$ in the compact $\omb$, there is at most
one $\omb$-valued solution to the equation $(M+D)$.
\end{stheorem}

{ \bf Remark : } Theorem \ref{unicit3} is optimal in the following
sense~: if $\omb$ is a geodesic ball of radius $r \ge \pi / (2 \sqrt K)$,
then it is possible to exhibit two different martingales (in particular,
two different solutions of the equation $(M+D)$ with $f=0$) which have
the same terminal value. \\
Let us consider the classical example of the sphere $\S^2$. We embed
$\S^2$ in $\R^3$, as the sphere of radius $1$. We call $N=(0,0,1)$ the
northern pole and we take $\omb = {\cal B}_{\rho_0}$, the geodesic ball
centered at $N$ of radius $\pi /(2 \sqrt K) = \pi /2$. Note
that $\omb$ is nothing but the northern hemisphere, containing the equator
$$E= \{ (x,y,z) / x^2+y^2=1 \ \hbox{ and } \ z=0 \}.$$
Now let $(B_t)$ be a BM starting at $(1,0,0)$, moving along the equator
$E$ and stopped when it reaches the plane $\{ x=0 \}$. Define another BM 
$(B'_t)$ moving along $E$ by reflecting $(B_t)$ with respect to the plane
$\{ x=0 \}$ (in particular, it starts from $(-1,0,0)$ and is stopped when
it reaches the plane $\{ x=0 \}$). These two processes are martingales
on the sphere, since they are BM moving along a geodesic; moreover, it is
obvious that they have the same terminal value (the point $(0,1,0)$ or
$(0,-1,0)$). Hence, in this case, the uniqueness doesn't hold. \\
Note that this result can be extended to manifolds which have closed 
geodesics (see for instance Proposition 2.2.2 of \cite{mpl}).

\bigskip
We conclude the uniqueness part by giving a consequence of the calculus
achieved in the two preceding subsections, which will be useful in the 
proof of existence.
\begin{sprop}
\label{p-integ}
In the two cases above (general connection and regular geodesic balls), 
there is a $q_0>1$ such that the submartingale $(S_t)_t$ is in 
$L^q(\Omega)$ for $1<q<q_0$.
\end{sprop}

\section{Existence}
\label{par14}

In this section we are given an $\omb$-valued random variable $U$
and we want to construct a pair of processes $(X,Z)$, satisfying
equation $(M+D)$, with $X$ in $\omb$ and terminal value $U$.
We limit ourselves to the case of a Wiener probability space.
{\bf This part is so similar as the existence part of
\cite{blache03} that we just give the main results and the changes
to make in the proofs.}

We recall the strategy of the proof in \cite{blache03}~:

\noindent
1. Simplify the problem by considering only terminal values which
can be expressed as functions of the diffusion $B^y$ at time $T$,
i.e. $U=F(B_T^y)$. This step needs to pass
through the limit in equation $(M+D)$.

\noindent
2. Solve a Pardoux-Peng BSDE with parameter to construct a pair of
processes in $\R^n \times \R^{nd_w}$ which is close to being a solution
of $(M+D)$ with $X_T=U$.

\noindent
3. Show that under an additional condition on $f$ the solution of
the preceding BSDE is a solution of the BSDE $(M+D)$ on a small time
interval.

\noindent
4. Use the convex function $\Psi$ to show that we have a solution of
$(M+D)$ on the whole time interval $[0;T]$.

\noi In fact, for technical reasons we suppose in the two last steps that
$f$ is sufficiently regular; then the proof of the existence is
completed with the last subsection~:

\noindent
5. Solve BSDE $(M+D)$ for general $f$ using classical approximation
methods.

\medskip
Note that we usually work within local coordinates in $\R^n$, i.e. we
consider that $\omb \subset O \subset \R^n$.

\subsection{Reduction of the problem}
\label{par14.1}
For an $\omb$-valued random variable $U^l$, we denote by $(X^l,Z^l)$
the solution of BSDE $(M+D)$ with terminal value $U^l$, and such that
$X^l$ remains in $\omb$.

\begin{sprop}
\label{generalU2}
If $U^l$ tends to $U$ in $L^2(\Omega)$ as $l \rightarrow \infty$, then
the processes $X^l$ tend to an $\omb$-valued process $X$ for the square 
norm $\esp ( \sup_{t \in [0,T]} \vert X^l_t - X_t \vert^2 )$, and the
processes $Z^l$ to a process $Z$ for the square norm
$$\esp \left( \int_0^T \Vert Z^l_t - Z_t \Vert^2 dt \right).$$
Moreover, $(X,Z)$ is the $\omb$-valued solution of BSDE $(M+D)$ with
terminal value $U$.
\end{sprop}

\begin{demo}
For integers $l,m$, we put
$$\tilde X^{l,m}=(X^l,X^m) \ \ \hbox{ and } \ \
  \tilde Z^{l,m} = \left(
                 \begin{array}{c}
                    Z^l   \\
	            Z^m
                  \end{array}
                \right).$$
We first deal with the existence of $X$. An application of Doob's $L^p$
inequality to the submartingale 
$(S_t)_t=(\exp (A_t) \Psi (\tilde X_t))_t$
$(S_t)_t$, with $1<q<q_0$ as in Proposition \ref{p-integ}, shows that, since 
$(U^l)_l$ is Cauchy, the sequence $(X^l)_l$ is Cauchy for the sup norm
$$ \delta^{(2)} \left( (X^l_t), (X_t) \right) =
   \sqrt {\esp \left( \sup_{t \in [0,T]} \delta^2 (X^l_t, X_t)
                                                 \right)}.$$
The proof is similar to the one of Lemma 4.1.2 in \cite{blache03}, but
we recall briefly the method for the convenience of the reader.
Since $\Psi \approx \delta^p$, we have
\begin{eqnarray}
\esp \left( \sup_t \delta^p(X^l_t,X^m_t) \right)
& \le & C \esp \left( \sup_t
        \left( e^{A_t} \Psi(X^l_t,X^m_t) \right)^q \right)^{\frac{1}{q}}\nonumber  \\
& \le & C \esp \left( e^{q A_T} \Psi^q(X^l_T,X^m_T) \right)^{\frac{1}{q}} \nonumber  \\
& \le & C \esp \left( e^{q \tilde q A_T} \right) ^\frac{1}{q \tilde q}
         \esp \left( \delta^{\frac{q \tilde q p}{\tilde q -1}}(U^l,U^m) \right) 
	    ^\frac{\tilde q -1}{q \tilde q}                             \nonumber   \\
\label{major2}
\end{eqnarray}
The constant $C$ above is allowed to vary from one inequality to
another, but it depends only on $T$, $\omb$ and $\Psi$ (but not on
the processes $(X^l)$).
The second inequality is Doob's $L^q$ one applied to the submartingale
$(\exp(A_t) \Psi(X_t,X'_t))_t;$ and the third one is H\" older's one
with $\tilde q$.\\
We choose $q$ and $\tilde q$ such that $q, \tilde q>1$ and $q \tilde q < q_0$
with $q_0$ as in Proposition \ref{p-integ}, so that $\exp(q \tilde q A_T)$ is
integrable. Note that $\esp (\exp(q \tilde q A_T))$ is uniformly bounded by
Lemmas \ref{2integ} and \ref{leminteg}.
Moreover, there are positive constants $\tilde c$, $\gamma_1$ and $\gamma_2$ such 
that, for any $\gamma >1$ and any variables $X_1$ and $X_2$ in $\omb$,
$$\tilde c \esp ( \de^\gamma(X_1,X_2))^\frac{1}{\gamma_1} \le 
\esp ( \de^2(X_1,X_2))
\le  \frac{1}{\tilde c} \esp ( \de^\gamma(X_1,X_2))^\frac{1}{\gamma_2};$$
this easily results from the boundedness of $\delta$ or H\" older's inequality.
At the end we get, for a constant $\zeta >0$ 
$$ \delta^{(2)} \left( (X^l_t), (X_t) \right) \le 
   C \esp \left( \delta^{2}(U^l,U^m) \right)^{\zeta}.$$
Now if $(U^l)_l$ is Cauchy in $L^2$, the sequence $(X^l)_l$ is clearly Cauchy for 
the sup norm. Thus $(X^l)_l$ converges to a process $X$.

\bigskip
Then we seek a process $Z$.
We consider here the function $\hat \Psi$ defined by~:
$$
\left\{
\begin{array}{ll}
\hat \Psi = \Psi \approx \de^2 
&  \hbox{ if a general connection is used (see \eqref{defpsiemery})}; \\
\hat \Psi = \de^2 
&  \hbox{ in the case of the Levi-Civita connection.} 
\end{array}
\right.
$$
For this function $\hat \Psi$, there is $a>0$ and $b>0$ such that
$$\forall (x,x') \in \omb \times \omb, \ \forall (z,z')
   \in T_xM \times T_{x'}M,$$
\begin{equation}
{}^t
 \left(
   \begin{array}{c}
     z  \\
     z'
   \end{array}
 \right)
\Hess \hat \Psi (x,x')
 \left(
   \begin{array}{c}
     z  \\
     z'
   \end{array}
 \right)
\ge \al \left\vert \overset{x'}{\underset{x}{\Vert}} z - z' \right\vert_r^2
    - \beta \hat \Psi(x,x') (\vert z \vert_r^2 + \vert z' \vert_r^2).
\label{minhesspsi2}
\end{equation}
In the case of a general connection, this is a consequence of
\eqref{minhesspsi}.
If the Levi-Civita connection is used, \eqref{minhesspsi2} results from
\eqref{2der1} and \eqref{derkpos}, remarking that, for
$u \in T_xM \times T_{x'}M$,
$$ \Hess \delta^2 (x,x')<u,u> = 2 \left( \delta(x,x')
     \Hess\delta (x,x')<u,u> + (\delta'(x,x')<u>)^2 \right).$$
Now apply It\^o's formula \eqref{2ito1} to $\hat \Psi(\tilde X^{l,m})$ and 
use estimates \eqref{minhesspsi2} and \eqref{2majdpsi} with $\ep=\al / 2$ 
to write
\begin{eqnarray*}
\frac{\al}{2} \esp \int_0^T \left\Vert
  \overset{X_s^l}{\underset{X_s^m}{\Vert}} Z_s^m - Z_s^l
                                         \right\Vert_r^2 ds
&\le& \esp \left( \hat \Psi(\tilde X^{l,m}_T) + \hat \Psi(\tilde X^{l,m}_0)
       \right)                                                       \\
&   & + C_\ep \esp \int_0^T \hat \Psi(\tilde X_s^{l,m})
           ( 1 + \Vert Z_s^l \Vert_r + \Vert Z_s^m \Vert_r) ds        \\
&   & + \beta \esp \int_0^T \hat \Psi(\tilde X_s^{l,m})
             (\Vert Z_s^l \Vert_r^2 + \Vert Z_s^m \Vert_r^2) ds.
\end{eqnarray*}
Since $(X_l)_l$ is Cauchy for the sup $L^2$-norm, the three expectations
on the right tend to zero as $l, m \rightarrow \infty$. Indeed, for the 
first expectation, it suffices to note that $\hat \Psi \approx \de^2$, and
for the last two, to use H\" older's inequality and the fact that
$$\int_0^T (\Vert Z_s \Vert_r^2 + \Vert Z'_s \Vert_r^2) ds$$
has some uniformly bounded exponential moments (from
Lemmas \ref{2integ} and \ref{leminteg}).\\
An application of the inequality \eqref{2tp2} and an argument similar to that
used just above show that, for the Euclidean norm, 
$$
 \esp \int_0^T \Vert Z_s^m - Z_s^l \Vert ^2 ds
    \overset{l,m \rightarrow \infty}{\longrightarrow} 0.
$$
Then $(Z^l)_l$ is a Cauchy sequence for the square
norm $\esp ( \int_0^T \Vert Z^l_t - Z^m_t \Vert^2 dt)$ and by completeness,
there is a limit process $Z$.

To complete the proof, it remains, by passing through the limit, to show that
$(X,Z)$ is the solution of BSDE $(M+D)$ with $X_T=U$. This is easy and has yet been
completed in the Second Step of the proof of Proposition 4.1.4 in
\cite{blache03}.
\end{demo}

Using the density of the space of all functionals
$$
\begin{array}{rl}
\{ G(W_{t_1}, W_{t_2}, \ldots , W_{t_q}), 
&  0<t_1 < \ldots < t_q \le T ;    \\
&   G : \R^{q d_W} \rightarrow \omb \ \hbox{ smooth, constant off a compact set } \}
\end{array}
$$
in $L^2(\cF_T ; \omb)$ and an argument of conditioning (see the end of Section 4.1
in \cite{blache03}), it turns out that it
suffices to solve equation $(M+D)$ with a terminal
value $U$ that can be written $F(B_T^y)$ for smooth 
$F : \R^d \rightarrow \omb$, constant off a compact set.

\subsection{The solution for terminal values $U=F(B_T^y)$}
\label{par14.2}
The proof of the existence of such a solution corresponds to the 
steps 2, 3 and 4 given at the beginning of Section \ref{par14}. In
\cite{blache03}, they are dealt with in Subsections 4.2, 4.3, 4.4, 4.5. An 
accurate examination of these subsections shows that the
proofs go essentially the same for the two cases here.
The only change to make is in the proof of Proposition 4.5.1, when we apply
Doob's inequality to the submartingale $(S_t)_t$~: instead of an
$L^2$ inequality, we have to consider an $L^p$ inequality for $p$ near $1$,
as in Proposition \ref{p-integ}. 
Therefore, under the additional condition
$$(H_s) \ \ f \hbox{ is pointing strictly outward on the boundary }
\partial \omb \ \hbox{ of } \omb,$$
which means that (for the Riemannian inner product $( \cdot \vert \cdot)_r$)
\begin{equation}
 \forall (b,x,z) \ : \ x \in \partial \omb, \ \
   \underset{b,x,z}{\inf}(D \chi(x) | f(b,x,z))_r \ge \zeta >0,
\label{Hs2}
\end{equation}
we have the following existence result.
\begin{sprop}
\label{existence2s}
Consider the BSDE $(M+D)$
with a terminal value $U = F(B_T^y)$ in $\omb= \{ \chi \le c \}$. Suppose 
that $f$ is a $C^3$ function with bounded derivatives, which verifies 
conditions \eqref{lip11}, \eqref{upperboundf11} and $(H_s)$; suppose 
moreover that $\chi$ is strictly convex (i.e. $\Hess \chi$ is positive 
definite). Then
\begin{item}
(i) If a general connection is used, there is a solution $(X,Z)$ to
$(M+D)$ with $X \in \omb$.
\end{item}
\begin{item}
(ii) If $\omb \subset {\cal B}_\rho$  (a regular geodesic ball of radius 
$\rho$ and center $o$), there is again a solution $(X,Z)$ to $(M+D)$ with 
$X \in \omb$.
\end{item}
\end{sprop}

{ \bf Remark :} 
If $\omb = {\cal B}_\rho$ in $(ii)$, the condition on $\chi$ is fulfilled. 
Indeed, ${\cal B}_\rho =\{ \chi \le c \}$ for a $c$ such that 
$0 \le c \le \pi^2 / (4K)$ and the smooth convex function (on 
${\cal B}_\rho$) $\chi(x) = \de^2(o,x)$. The strict convexity of $\chi$
on ${\cal B}_\rho$ comes from \eqref{2der1}, \eqref{minhessdelta} and the
classical formula, for $z \in T_xM$,
$$ \Hess \chi (o,x)<z,z> = 2 \left( \delta(o,x)
     \Hess\delta (o,x)<z,z> + (\delta'(o,x)<z>)^2 \right):$$
these formulas imply that, for $z=v+w \in T_xM$ (where $v$ and $w$ are 
the tangential and orthogonal components of $z$ (see Subsection 
\ref{2estimdist})), 
\begin{eqnarray*}
 \Hess \chi (o,x) <z,z> 
&\ge& \sqrt K \de(o,x) \cot (\sqrt K \de(o,x)) 
        \vert w \vert^2_r + \vert v \vert^2_r            \\
&\ge& \eta \vert z \vert_r^2 \ \hbox{ on ${\cal B}_\rho$, for $\eta>0$.}
\end{eqnarray*}
In this case, the condition \eqref{Hs2} is equivalent to
$$  \forall (b,x,z) \ : \ x \in \partial \omb, \ \
    \underset{b,x,z}{\inf}  \ \de'(o,x) < f(b,x,z)> \ge \zeta >0.$$

\bigskip
To get the existence in a general framework, it remains to extend  
Proposition \ref{existence2s} to general terminal values $U$
and nonregular $f$, verifying the weaker condition
$$(H) \ \ f \hbox{ is pointing outward on the boundary of } \omb.$$
Now it means that
$
 \forall (b,x,z) \ : \ x \in \partial \omb, \ \
   (D \chi(x) | f(b,x,z))_r \ge 0.
$
Note that this is a natural condition since it is necessary in the 
deterministic case (i.e. when the terminal value $U$ is deterministic 
and $Z=0$); moreover, in the case of a regular geodesic ball, it is
equivalent to require that $ \de'(o,x) < f(b,x,z)>$ be only 
nonnegative.
This generalization follows from Proposition \ref{generalU2}, and 
Subsection 4.6 of \cite{blache03}. At the end, we get the existence 
result of Theorem \ref{existence12}.

\section{Applications}
\label{par15}

\subsection{The martingale case}
\label{par15.1}
In the case of a vanishing drift $f$, solving equation $(M+D)$ is
equivalent to finding a martingale on $M$ with terminal value $U$. Then
we recover Kendall's results for regular geodesic balls, stated in
\cite{kend90}.

\subsection{Case of a random terminal time}
\label{par15.2}
This case is the same as in \cite{blache03}. So we just give the results
obtained. The equation which we now study is
$$(M+D)_\tau
\left\{
      \begin{array}{l}
          d X_t = Z_t d W_t + \left( - \half \Gamma_{jk}(X_t)
            ([Z_t]^k \vert [Z_t]^j) + f(B_t^y,X_t,Z_t) \right) d t  \\
	  X_\tau=U^\tau                                         \\
      \end{array}
    \right.   $$
where
$U^\tau$ is a $\omb$-valued, ${\cal F}_\tau$-measurable random variable
and $\tau$ is a stopping time which verifies the exponential 
integrability condition
\begin{equation}
  \exists \xi > 0 :  \esp (e^{\xi \tau})<\infty.
\label{2integst}
\end{equation}
In this case, we need an additional restriction on $f$ to keep the
integrability condition 
\begin{equation}
 \esp \left( e^{\lambda \tau + \mu \int_0^\tau
  (\Vert Z_s \Vert_r^2 + \Vert Z'_s \Vert_r^2) ds} \right) < \infty.
\label{integstS_t}
\end{equation}
The condition is that $f$ should be "small", in the following sense~:
there is an $h<\xi$ such that if the Lipschitz constant $L$ and the 
bound $L_2$ (respectively in \eqref{lip11} and \eqref{upperboundf11}) of 
$f$ are smaller than $h$, then the integrability condition in 
\eqref{integstS_t} holds (in particular, it implies that $\la < \xi$).
In this case, $(S_t)_{0 \le t \le \tau}$ remains a submartingale, and 
the consequence is the
\begin{stheorem}
\label{2existence4}
We consider the BSDE $(M+D)_\tau$ with $\tau$ a stopping time verifying the
integrability condition \eqref{2integst}; the function $\chi$ used to
define the domain $\omb$ is supposed as usual to be strictly convex.
Then if $f$ verifies conditions \eqref{lip11}, \eqref{upperboundf11}, $(H)$
and moreover is "small" (in the sense defined above), this BSDE has a
unique solution $(X,Z)$, in the same cases as in Theorem \ref{existence12}.
\end{stheorem}
Note that if the stopping time $\tau$ is bounded a.s., then 
\eqref{integstS_t} holds (and therefore the existence and uniqueness of a
solution) without supposing that $f$ should be "small".

We end this paper by recalling briefly from \cite{blache03} the applications
to PDEs. The reader is referred to paragraphs 5.4 and 5.5 of
\cite{blache03} for further details.

\subsection{Application to nonlinear PDEs}
\label{par15.4}
Under some conditions on the coefficients $\si$ and $b$ in the definition
\eqref{sde2} of the diffusion $B^x$, we can think of this diffusion as a
Brownian Motion on a Riemannian manifold $(N,g)$.
Let $\overline M_1$ be a compact submanifold of $N$, with boundary
$\partial M_1$ and interior $M_1$. 
Given a regular mapping
$\overline \phi : \partial M_1 \rightarrow \omb \subset M,$
we wish to find a mapping $\phi : \overline M_1 \rightarrow \omb$
which solves the Dirichlet problem
$$
(D)
\left\{
\begin{array}{ccc}
{\cal L}_M \phi(x) - f(x, \phi(x), \nabla \phi(x) \si(x))=0
& , & x \in M_1                                        \\
\phi(x) = \overline \phi (x)
& , & x \in \partial M_1
\end{array}
\right.
$$
where ${\cal L}_M \phi$ is the tension field of the mapping $\phi$.

Let $\zeta$ denote the first time $B^x$ hits the boundary; we assume that 
$\zeta$ verifies an integrability condition like \eqref{2integst}.
Using the same Wiener process $W$ with which we constructed $B^x$, we can
solve according to Theorem \ref{2existence4} the BSDE $(M+D)_\zeta$ with
terminal value $\overline \phi (B^x_\zeta)$. Let
$(X_t^x, Z_t^x)_{0 \le t \le \zeta}$ be the unique solution and put
$\phi(x):=X_0^x$. Then
under sufficient regularity on $\phi$, it is not difficult to verify that
$\phi$ is a solution to the Dirichlet problem $(D)$. 

Note also that when $f(b,x,z)=f(b,x)$ and is written as
$f(b,x)  = D_2 G(b,x)$ (the differential of $G$ with respect to the
second variable), the elliptic nonlinear PDE in the Dirichlet
problem $(D)$ is associated with a variational problem; but it seems harder
to associate variational problems for more general $f$.

\bigskip
We conclude by remarking that we can solve, as in \cite{blache03},
the heat equation  associated with the elliptic problem $(D)$.

\bibliography{biblio}
\bibliographystyle{abbrv}

\end{document}